\documentclass[12pt,a4paper, oneside,reqno]{article}
\usepackage[utf8]{inputenc}
\usepackage[english]{babel}
\usepackage[left=1in,right=1in,top=1in,bottom=1in]{geometry}
\usepackage{amsfonts,amssymb,amsmath,amsthm}
\usepackage{verbatim}

\usepackage[symbol]{footmisc} 
\usepackage{longtable}
\usepackage{enumerate}
\usepackage{hyperref}
\usepackage{tkz-graph,comment}

\usepackage{fancyhdr} 
\theoremstyle{plain}
\newtheorem{theorem}{Theorem}
\newtheorem{lemma}[theorem]{Lemma}
\newtheorem{proposition}[theorem]{Proposition}

\newtheorem{definition}[theorem]{Definition}

\usetikzlibrary{arrows}

\newtheorem*{theoremA1}{Theorem A1}
\newtheorem*{theoremA2}{Theorem A2}
\newtheorem*{theoremA3}{Theorem A3}
\newtheorem*{theoremA4}{Theorem A4}
\newtheorem*{theoremG1}{Theorem G1}
\newtheorem*{theoremG2}{Theorem G2}
\newtheorem*{theoremG3}{Theorem G3}
\newtheorem*{theoremG4}{Theorem G4}

\usepackage{fouriernc} 

\usepackage{cite}

\begin{document}


 \bigskip

\noindent{\Large
The algebraic and geometric classification of \\
derived Jordan and bicommutative algebras}
\footnote{The  authors would like to thank the SRMC 
$\big($Sino-Russian Mathematics Center in Peking University, Beijing, China$\big)$ for its hospitality and excellent working conditions, where some parts of this work have been done.
The work is supported by 
FCT  2023.08031.CEECIND and UID/00212/2025.}

 \bigskip

\begin{center}

 {\bf
Hani Abdelwahab\footnote{Department of Mathematics, 
 Mansoura University,  Mansoura, Egypt; \ haniamar1985@gmail.com}, 
   Ivan Kaygorodov\footnote{CMA-UBI, University of  Beira Interior, Covilh\~{a}, Portugal; \    kaygorodov.ivan@gmail.com}   \&
   Roman Lubkov\footnote{Department of Mathematics and Computer Science, Saint Petersburg State University, Russia; r.lubkov@spbu.ru, romanlubkov@yandex.ru}  
}

\end{center}

\

\noindent {\bf Abstract:}
{\it  
We developed a new proper method for classifying $n$-dimensional derived Jordan algebras, 
and apply it to the classification of $3$-dimensional derived Jordan algebras. 
As a byproduct, we have the algebraic classification of $3$-dimensional metabelian commutative algebras and $3$-dimensional  derived commutative associative algebras. 
After that, we introduced a method of classifying $n$-dimensional bicommutative algebras, based on the classification of $n$-dimensional derived commutative associative algebras, and applied it to the classification of $3$-dimensional bicommutative algebras.
The second part of the paper is dedicated to the geometric classification of $3$-dimensional  metabelian commutative, derived commutative associative, derived Jordan and bicommutative algebras. 
}

 \bigskip 

\noindent {\bf Keywords}:
{\it 
Jordan algebras, 
bicommutative   algebras,
metabelian algebras,
algebraic classification,
geometric classification.}

\bigskip 

 \
 
\noindent {\bf MSC2020}:  
17A30 (primary);
17D25,
14L30 (secondary).

	 \bigskip

\ 


\tableofcontents 

\section*{Introduction}
The algebraic classification $\big($up to isomorphism$\big)$ of algebras of dimension $n$ of a certain variety
defined by a family of polynomial identities is a classic problem in the theory of non-associative algebras, see \cite{k23, jp3, jp2, akk,aak24,afm,kkp,kpv22,kpv20,AKL}.
There are many results related to the algebraic classification of small-dimensional algebras in different varieties of
associative and non-associative algebras.
For example, algebraic classifications of 
$2$-dimensional algebras,
$3$-dimensional Poisson algebras,
$4$-dimensional alternative  algebras,
$5$-dimensional symmetric Leibniz algebras 
and so on have been given.
 Deformations and geometric properties of a variety of algebras defined by a family of polynomial identities have been an object of study since the 1970s, see \cite{ben,BC99,FKS,GRH,FKS25, jp1,jp2,kz}  and references in \cite{k23,l24,MS}.

\medskip 

Jordan algebras are the most popular type of non-associative algebras\footnote{Due to many connections between associative and Lie algebras, we consider both as staying outside of the "non-associative world".}.
The first Jordan algebras appeared in a classical paper by  Jordan, von Neumann, and  Wigner in 1934 \cite{JVW}, and until now they have been under active investigation.
The most famous generalizations of Jordan algebras are the following:
non-commutative Jordan algebras $\big($see, \cite{aak24} and references therein$\big)$;
terminal algebras $\big($see, \cite{kkp} and references therein$\big)$;
commutative $\mathfrak{CD}$-algebras $\big($see, \cite{jkk} and references therein$\big)$;
axial algebras  $\big($see, \cite{BG} and references therein$\big)$;
structurable  algebras $\big($see, \cite{DM} and references therein$\big);$ and so on. 
In the present paper, we introduce another generalization of Jordan algebras. 
Namely, we are studying commutative derived Jordan algebras.
The idea of considering derived algebras arose from the study of the  lower central series in groups. 
Namely, for  a variety of algebras defined by a family of polynomial identities $\Omega$, 
we say that an algebra ${\rm A}$ is a derived $\Omega$-algebra if ${\rm A}^2$ is an $\Omega$-algebra.
One of the basic examples of derived $\Omega$-algebras is metabelian
$\big($also known as $2$-step solvable$\big)$ algebras, 
i.e., algebras satisfying the identity 
    $(xy)(zt)\ = \ 0.$
They are in the intersection of all non-(anti)commutative derived $\Omega$-algebras.
Metabelian algebras from a certain variety of algebras are on a certain interest now \cite{DAS,ms,kg,mv,p81,BCZ,SS}.
The notion of derived algebras plays an important role in many varieties of algebras. 
So, each derived algebra of a solvable Lie or Leibniz algebra will be nilpotent;
each derived algebra of a bicommutative 
$\big($or anti-bicommutative$\big)$ algebra $\big($under the Jordan product$\big)$ is commutative associative  \cite{DT,BBK}.
Hence, our results are given in subsection \ref{derasscom}  
$\big($classification of derived commutative associative algebras$\big)$ 
play an important role in the future classification of $3$-dimensional bicommutative algebras, given in subsection \ref{alg_bi}. 

 \medskip

The systematic study of bicommutative   algebras from the algebraic point of view started after a paper by Dzhumadildaev and Tulenbaev \cite{DT},  where they proved that 
bicommutative algebras are Lie admissible and derived commutative associative   under the Jordan product.
The main results from the theory of bicommutative algebras were obtained by 
Dzhumadildaev and Drensky with various co-authors. 
Namely, Dzhumadildaev,  Ismailov, and   Tulenbaev  constructed  bases for free bicommutative algebras and proved that the bicommutative operad is not Koszul \cite{DIT}. 
Dzhumadildaev and   Ismailov described all identities of commutator and anticommutator algebras of bicommutative algebras \cite{DI}.
Drensky found the bases of  polynomial identities of $2$-dimensional non-associative bicommutative algebras and showed that one of these algebras generates the whole variety of bicommutative algebras \cite{D19}. Later, he found analogies of classical results of the invariant theory of finite groups acting on polynomial algebras, namely analogies of the Endlichkeitssatz of Emmy Noether, the Molien formula, and the Chevalley-Shephard-Todd theorem. Finally, he describes the symmetric polynomials of free bicommutative algebras \cite{D23}.
Invariants of free metabelian bicommutative algebras were also studied in a paper by Öğüşlü and Fındık \cite{OF}.  
Drensky and  Zhakhayev   obtained  a positive solution to the Specht problem for  bicommutative algebras \cite{DZ}.
Finally, Drensky, Ismailov, Mustafa, and Zhakhayev established various superanalogs of classical theorems for bicommutative superalgebras in \cite{dimz}.
Ismailov, Mashurov, and Sartayev proved that 
every metabelian Lie algebra can be embedded into a bicommutative algebra with respect to the commutator product \cite{ims}.
Products of commutator ideals of bicommutative  algebras were studied in a paper by Kaygorodov, Mashurov, Nam, and  Zhang \cite{KMTZ}. 
Bai, Chen, and  Zhang proved that 
the Gelfand-Kirillov dimension of a finitely generated bicommutative algebra is a nonnegative integer \cite{bcz}.
Shestakov and Zhang gave an example of wild automorphisms of the free two-generated bicommutative algebra in \cite{SZ}.
Bialgebra theory for bicommutative algebras was established in \cite{BBR}. 
Towers introduced the Frattini theory for bicommutative algebras \cite{T25}.
Kolesnikov and Sartayev proved that 
the Hadamard product of the bicommutative  operad with the Novikov operad
coincides with their white Manin product \cite{ks25}.
The bicommutative operad also satisfies the Dong property \cite{DS}.
Bicommutative algebras are one of the most important subvarieties in Lie-admissible algebras $\big($about Lie-admissible algebras see \cite{DAS}$\big).$ 
Bicommutative algebras are algebras of biderivation-type \cite{BO}.
Burde and  Dekimpe described bicommutative algebraic structures on certain Lie algebras in \cite{BK,DKS}.
Some results regarding algebraic and  geometric classifications of small-dimensional nilpotent bicommutative algebras were obtained by 
Abdurasulov,     
Kaygorodov, Khudoyberdiyev, Páez-Guillán, and Voronin \cite{kpv20,kpv22,akk}.
We also mention that bicommutative algebras are a particular case of $\delta$-Novikov algebras for $\delta=0,$
and the present paper will finish the classification of $3$-dimensional $\delta$-Novikov algebras
$\big($for the case of $\delta \notin \big\{ 0,1\big\}$ see \cite{AKL}$\big).$

 \medskip

The present paper aims to give the complete algebraic and geometric classification of $3$-dimensional  algebras in the following varieties:
metabelian commutative, 
derived commutative associative, 
derived Jordan, 
and bicommutative algebras.
Namely, 
subsection \ref{met} is dedicated to the algebraic classification of $3$-dimensional metabelian commutative algebras, which is presented in Theorem A1; 
subsection \ref{derasscom} is dedicated to the algebraic classification of $3$-dimensional derived commutative associative algebras, which is presented in Theorem A2; 
subsection \ref{derJ} is dedicated to the algebraic classification of $3$-dimensional derived Jordan algebras,  which is presented in Theorem A3. The next subsection \ref{alg_bi} contains the algebraic classification of $3$-dimensional bicommutative algebras
$\big($Theorem A4$\big).$
Based on the obtained results, we present the geometric classifications: 

\medskip 
\noindent 
{\bf Theorem G1} states that the variety of $3$-dimensional metabelian  commutative algebras is $8$-dimensional with    
$2$ irreducible components and with $1$ rigid algebra.

\medskip 
\noindent 
{\bf Theorem G2} states that the variety of $3$-dimensional derived commutative associative algebras is $12$-dimensional with    
$2$ irreducible components and with $1$ rigid algebra. 

\medskip 
\noindent 
{\bf Theorem G3} states that the variety of $3$-dimensional derived Jordan algebras is $12$-dimensional with    
$7$ irreducible components and with $5$ rigid algebras. 

\medskip 
\noindent 
{\bf Theorem G4} states that the variety of $3$-dimensional bicommutative algebras is $10$-dimensional with $4$ irreducible components and with $1$ rigid algebra.


\section{The algebraic classification of  algebras}

All the algebras below will be over $\mathbb C$ and all the linear maps will be $\mathbb C$-linear.
For simplicity, every time we write the multiplication table of an algebra, the products of basic elements whose values are zero or can be recovered from the commutativity  or from the anticommutativity are omitted.
The notion of a nontrivial algebra means that the multiplication is nonzero.

The following propositions are well known and, for example, the multiplication tables of algebras can be found in \cite{aak24}.

\begin{proposition}[see \cite{aak24}]
Let ${\mathfrak J}$ be a nontrivial complex $2$-dimensional Jordan algebra.
Then ${\mathfrak J}$\ is isomorphic to one of the following algebras:

\begin{longtable}{lllllll}
${\mathfrak J}_{01}$ & $:$ & $e_{1}e_{1} = e_{1}$ & $e_{2}e_{2} = e_{2}$ & & & \\ 
${\mathfrak J}_{02}$ & $:$ & $e_{1}e_{1} = e_{1}$ & & & & \\ 
${\mathfrak J}_{03}$ & $:$ & $e_{1}e_{1} = e_{1}$ & $e_{1}e_{2} = e_{2}$ & & & \\ 
${\mathfrak J}_{04}$ & $:$ & $e_{1}e_{1} = e_{1}$ & $e_{1}e_{2} = \frac{1}{2}e_{2}$ & & & \\ 
${\mathfrak J}_{05}$ & $:$ & $e_{1}e_{1} = e_{2}$ & & & & \\ 
\end{longtable}\noindent
All algebras, excepting ${\mathfrak J}_{04},$ are associative. 
\end{proposition}

\begin{proposition}[see \cite{aak24}]
\label{3-dim assoc}Let ${\rm J}$ be a complex $3$-dimensional
commutative associative algebra. Then ${\rm J}$ is isomorphic to one of
the following algebras:

\begin{longtable}{llllll}
${\rm J}_{01}$ & $:$ & $e_{1}e_{1} = e_{1}$ & $e_{2}e_{2} = e_{2}$ & \\ 
${\rm J}_{02}$ & $:$ & $e_{1}e_{1} = e_{1}$ & $e_{1}e_{2} = e_{2}$ & \\ 
${\rm J}_{03}$ & $:$ & $e_{1}e_{1} = e_{1}$ & & \\ 
${\rm J}_{04}$ & $:$ & $e_{1}e_{1} = e_{2}$ & & \\ 
${\rm J}_{05}$ & $:$ & $e_{1}e_{2} = e_{3}$ & & \\ 
${\rm J}_{06}$ & $:$ & $e_{1}e_{1} = e_{2}$ & $e_{1}e_{2} = e_{3}$ & \\ 
${\rm J}_{07}$ & $:$ & $e_{1}e_{1} = e_{1}$ & $e_{2}e_{2} = e_{2}$ & $e_{3}e_{3} = e_{3}$ \\ 
${\rm J}_{08}$ & $:$ & $e_{1}e_{1} = e_{1}$ & $e_{2}e_{2} = e_{2}$ & $e_{2}e_{3} = e_{3}$ \\ 
${\rm J}_{09}$ & $:$ & $e_{1}e_{1} = e_{1}$ & $e_{1}e_{2} = e_{2}$ & $e_{1}e_{3} = e_{3}$ \\ 
${\rm J}_{10}$ & $:$ & $e_{1}e_{1} = e_{1}$ & $e_{1}e_{2} = e_{2}$ & $e_{1}e_{3} = e_{3}$ & $e_{2}e_{2} = e_{3}$ \\ 
${\rm J}_{11}$ & $:$ & $e_{1}e_{1} = e_{1}$ & $e_{2}e_{2} = e_{3}$ & \\ 
\end{longtable}
 
\end{proposition}

\begin{proposition}[see \cite{aak24}]
\label{3-dim Jord}Let ${\rm J}$ be a complex $3$-dimensional Jordan
algebra. Then ${\rm J}$ is a commutative associative algebra listed in
{\rm Proposition~\ref{3-dim assoc}} or isomorphic to one of the following algebras:

\begin{longtable}{l l l l l l l l}
${\rm J}_{12}$ & $:$ & $e_{1}e_{1} = e_{1}$ & $e_{2}e_{2} = e_{2}$ & $e_{3}e_{3} = e_{1} + e_{2}$ & $e_{1}e_{3} = \frac{1}{2}e_{3}$ & $e_{2}e_{3} = \frac{1}{2}e_{3}$ & \\ 
${\rm J}_{13}$ & $:$ & $e_{1}e_{1} = e_{1}$ & $e_{1}e_{2} = \frac{1}{2}e_{2}$ & $e_{1}e_{3} = e_{3}$ & & & \\ 
${\rm J}_{14}$ & $:$ & $e_{1}e_{1} = e_{1}$ & $e_{1}e_{2} = \frac{1}{2}e_{2}$ & $e_{1}e_{3} = \frac{1}{2}e_{3}$ & & & \\ 
${\rm J}_{15}$ & $:$ & $e_{1}e_{1} = e_{1}$ & $e_{2}e_{2} = e_{3}$ & $e_{1}e_{2} = \frac{1}{2}e_{2}$ & & & \\ 
${\rm J}_{16}$ & $:$ & $e_{1}e_{1} = e_{1}$ & $e_{2}e_{2} = e_{3}$ & $e_{1}e_{2} = \frac{1}{2}e_{2}$ & $e_{1}e_{3} = e_{3}$ & & \\ 
${\rm J}_{17}$ & $:$ & $e_{1}e_{1} = e_{1}$ & $e_{2}e_{2} = e_{2}$ & $e_{1}e_{3} = \frac{1}{2}e_{3}$ & $e_{2}e_{3} = \frac{1}{2}e_{3}$ & & \\ 
${\rm J}_{18}$ & $:$ & $e_{1}e_{1} = e_{1}$ & $e_{1}e_{2} = \frac{1}{2}e_{2}$ & & & & \\ 
${\rm J}_{19}$ & $:$ & $e_{1}e_{1} = e_{1}$ & $e_{2}e_{2} = e_{2}$ & $e_{1}e_{3} = \frac{1}{2}e_{3}$ & & & \\ 
\end{longtable}

\end{proposition}

\subsection{Metabelian commutative algebras}\label{met}

\begin{definition}
An algebra ${\rm A}$ is called a metabelian    algebra if it satisfies 
    $(xy)(zt) \ = \ 0.$
\end{definition}

\begin{theoremA1}\label{metabelian}
Let ${\rm A}$ be a complex $3$-dimensional metabelian commutative    algebra.
Then ${\rm A}$
is  
isomorphic to one of the following algebras:

\begin{longtable}{lclllll}
 ${\rm M}_{01}
 $ & $:$ & $e_{1}e_{1} = e_{2}$ & & \\ 
 ${\rm M}_{02}
 $ & $:$ & $e_{1}e_{2} = e_{3}$ & & \\ 
 ${\rm M}_{03}
 $ & $:$ & $e_{1}e_{1} = e_{2}$ & $e_{1}e_{2} = e_{3}$ & \\ 
 ${\rm M}_{04}^{\alpha}$
 &$:$ & $e_{1}e_{3}=e_{1}$&$ e_{2}e_{3}=\alpha e_{2}$ \\
 ${\rm M}_{05}
 $&$:$ & $e_{1}e_{3}=e_{1}$&$ e_{3}e_{3}=e_{2}$ \\
 ${\rm M}_{06}
 $&$:$ & $e_{1}e_{3}=e_{1}+e_{2}$&$ e_{2}e_{3}=e_{2}$ \\
 ${\rm M}_{07}
 $&$:$&$e_{1}e_{3}=e_{1}$&$e_{2}e_{2}=e_{1}$
\end{longtable}
\noindent
All algebras, excepting  ${\rm M}_{04}^{\alpha }\cong {\rm M}_{04}^{\alpha  ^{-1}},$
are non-isomorphic.\\ \noindent
Algebras ${\rm M}_{01},$ ${\rm M}_{02},$ and ${\rm M}_{03}$ are associative.

\end{theoremA1}

\begin{proof}
It is easy to see that we have only two opportunities: 
$\dim {\rm A}^{2}=2$ and $\dim {\rm A}^{2}=1.$ 
\medskip 

\noindent 
{\bf First}, we consider algebras with $\dim {\rm A}^{2}=2.$ 
The
multiplication table of ${\rm A}$ is defined as follows: 
\begin{longtable}{lclllll}
${\rm A}$&$
:$&$e_{1}e_{3}=C_{131}e_{1}+C_{132}e_{2}$ & 
$e_{2}e_{3}=C_{231}e_{1}+C_{232}e_{2}$ &
$e_{3}e_{3}=C_{331}e_{1}+C_{332}e_{2}$
\end{longtable}  \noindent If $\big(C_{131}, C_{232}\big) = (0,0)$ and $C_{132} C_{231} = 0$, then $\rm A$ is an associative algebra, giving rise to the algebra ${\rm M}_{03}$. 
\medskip

We now turn to the case where $\big(C_{131}, C_{232}\big) \neq (0,0)$ or $C_{132} C_{231} \neq 0$. 
Since $\dim {\rm A}^{2} = 2$, the matrix 
\[
\mathfrak{C} =
\begin{pmatrix}
C_{131} & C_{132} \\ 
C_{231} & C_{232} \\ 
C_{331} & C_{332}
\end{pmatrix}
\]
necessarily has rank $2$.
For convenience, we encode the multiplication by the triple
$\theta = (B_1, B_2, B_{3})$, where
\begin{longtable}{lcl}
$B_{1} $&$=$&$C_{131}\triangledown_{13}+C_{231}\triangledown_{23}
+ C_{331}\triangledown_{33},$ \\
$B_{2} $&$=$&$C_{132}\triangledown_{13}+C_{232}\triangledown_{23}
+ C_{332}\triangledown_{33},$\\
$B_{3} $&$=$&$0.$
\end{longtable}%

\noindent \ Here $\triangledown_{ij}$ denotes the symmetric bilinear form satisfying
\[
\triangledown_{ij}(e_p, e_q) = 
\begin{cases}
1 & \text{if } (p,q) = (i,j) \text{ or } (p,q) = (j,i), \\[2pt]
0 & \text{otherwise}.
\end{cases}
\] Then the product in $\mathrm{A}$ is given by
\[
x \cdot y = \theta(x,y) = B_1(x,y)\, e_1 + B_2(x,y)\, e_2, \quad x,y \in \mathrm{A}.
\]
Both $B_1$ and $B_2$ are nonzero because $\dim \mathrm{A}^2 = 2$.

\medskip\noindent 
The group $\mathrm{GL}_3(\mathbb{C})$ acts on algebra structures via a change of the basis. 
For $\phi \in \mathrm{GL}_3(\mathbb{C})$, the transformed structure is
\[
\theta \ast \phi = (B_1', B_2', B_3'), \qquad 
(\theta \ast \phi)(x, y) := \phi^{-1}\big(\theta(\phi(x), \phi(y))\big), \quad x, y \in \mathrm{A}.
\]
By construction, $\phi\colon (\mathrm{A}, \theta) \to (\mathrm{A}, \theta \ast \phi)$ is an algebra isomorphism, 
so $(\mathrm{A}, \theta) \cong (\mathrm{A}, \theta \ast \phi)$. 
Conversely, any isomorphism between such algebras arises from a change of basis; 
thus, classifying these algebras up to isomorphism amounts to describing the orbits of this action.

\medskip\noindent
To simplify the structure constants, we may choose a transformation of the form
\[
\phi =
\begin{pmatrix}
a_{11} & 0 & a_{13} \\ 
a_{21} & a_{22} & a_{23} \\ 
0 & 0 & a_{33}
\end{pmatrix}, \qquad a_{11} a_{22} a_{33} \neq 0.
\]
Then the transformed component $B_1'$ becomes
\begin{center}
$B_1' \ = \  
\frac{a_{33}(C_{131}a_{11}+C_{231}a_{21})}{a_{11}} \triangledown_{13}
+ \frac{C_{231}a_{22}a_{33}}{a_{11}} \triangledown_{23}
+ \frac{a_{33}(2C_{131}a_{13}+2C_{231}a_{23}+C_{331}a_{33})}{a_{11}} \triangledown_{33}.$\end{center}

\noindent From here, we may assume 
$B_{1}\in \big\{ \triangledown
_{13},\triangledown_{23},\triangledown_{33}\big\} $, i.e., 
\begin{center}
    $\big(
C_{131},C_{231},C_{331}\big) \in \big\{ \big( 1,0,0\big) ,\ 
\big(0,1,0\big) ,\ \big( 0,0,1\big) \big\} $.

\end{center}

\begin{itemize} 
\item $\big( C_{131},C_{231},C_{331}\big) =\big( 1,0,0\big) $.

\begin{itemize}
\item If $%
C_{232}\notin \big\{ 0,1 \big\}$, we choose $\phi $ to be the following matrix:%
\begin{equation*}
\begin{pmatrix}
1 & 0 & 0 \\ 
\frac{C_{132}}{1-C_{232}} & 1 & -\frac{C_{332}}{2C_{232}} \\ 
0 & 0 & 1%
\end{pmatrix}%
.
\end{equation*}%
Then $\theta \ast \phi =\big( \triangledown_{13},\ C_{232}\triangledown_{23},\ 0\big) $. So we get the algebras ${\rm M}_{04}^{\alpha \notin \{
0,1\}} $. Furthermore, ${\rm M}_{04}^{\alpha }\cong {\rm M}_{04}^{\beta}$ if 
$ \alpha =  \beta ^{-1}$.
 
    \item 
If $C_{232}=0$, then $C_{332}\neq 0$ since the matrix $%
\mathfrak{C}
$ has rank $2$. Choose $\phi $ as follows: 
\begin{equation*}
\begin{pmatrix}
1 & 0 & 0 \\ 
C_{132} & C_{332} & 0 \\ 
0 & 0 & 1%
\end{pmatrix}%
.
\end{equation*}%
Then $\theta \ast \phi =\big( \triangledown_{13},\ \triangledown_{33},\ 0\big) $. Thus we obtain the algebra ${\rm M}_{05}$. 

\item If $%
C_{232}=1$ and $C_{132}=0$, we choose $\phi $ to be the following matrix:%
\begin{equation*}
\begin{pmatrix}
1 & 0 & 0 \\ 
0 & 1 & -\frac{1}{2}C_{332} \\ 
0 & 0 & 1%
\end{pmatrix}%
.
\end{equation*}%
Then $\theta \ast \phi =\big( \triangledown_{13},\ \triangledown_{23},\ 0\big) $. So we get the algebra ${\rm M}_{04}^{1}$. 
\item If $%
C_{232}=1$ and $C_{132}\neq 0$, we choose $\phi $ to be the following matrix:%
\begin{equation*}
\begin{pmatrix}
1 & 0 & 0 \\ 
0 & C_{132} & -\frac{1}{2}C_{332} \\ 
0 & 0 & 1%
\end{pmatrix}%
.
\end{equation*}%
Then $\theta \ast \phi =\big( \triangledown_{13},\triangledown
_{13}+\triangledown_{23},0\big) $. Therefore we get the algebra $%
{\rm M}_{06}$.
\end{itemize}

\end{itemize}

\begin{itemize}
\item $\big( C_{131},C_{231},C_{331}\big) =\big( 0,1,0\big) $. Then $ C_{232}\neq 0$ or $C_{132} C_{231} \neq 0$.

\begin{itemize}

\item If $C_{132}\neq 0,$ 
we choose $\phi ^{\prime }$ to be the following matrix:%
\begin{equation*}
\begin{pmatrix}
{C^{-\frac{1}{2}}_{132}} & 0 & -\frac{1}{2} {C_{332}}{C_{132}^{-\frac{3}{2}}} \\ 
0 & 1 & 0 \\ 
0 & 0 &  {C^{-\frac{1}{2}}_{132}}%
\end{pmatrix}%
.
\end{equation*}%
Then $\theta \ast \phi ^{\prime }=\big( \triangledown_{23},\triangledown
_{13}+\alpha \triangledown_{23},0\big) $ where 
$\alpha =C^{-\frac{1}{2}}_{132}C_{232}$. Consider the following matrices and 
let $\phi ^{\prime\prime }=\varphi_1$ if $\alpha ^{2}+4=0$ 
and $\phi ^{\prime\prime }=\varphi_2$ if $\alpha ^{2}+4\neq 0$:%
\begin{equation*}
\varphi_1=\begin{pmatrix}
-\frac{2}{\alpha } & \frac{2}{\alpha } & 0 \\ 
0 & 1 & 0 \\ 
0 & 0 & \frac{2}{\alpha }%
\end{pmatrix}%
, \ 
\varphi_2=\begin{pmatrix}
\frac{1}{2} \big(\sqrt{\alpha ^{2}+4}-\alpha\big) & -\frac{1}{2} \big(\sqrt{\alpha
^{2}+4}+\alpha\big) & 0 \\ 
1 & 1 & 0 \\ 
0 & 0 & \frac{1}{2}\big(\sqrt{\alpha ^{2}+4}-\alpha\big)
\end{pmatrix}%
.\allowbreak
\end{equation*}%
Then \begin{center}
$\theta \ast \phi ^{\prime }\phi ^{\prime \prime }\in \big\{ \big(
\triangledown_{13},\triangledown_{13}+\triangledown_{23},0\big)
, \ \big( \triangledown_{13},\beta \triangledown_{23},0\big) \big\} $.\end{center}
Hence, there are no new algebras.

\item If $C_{132}=0,$ then $C_{232}C_{332}\neq 0$. Let $\phi$ be the following matrix:%
\begin{equation*}
\phi=\begin{pmatrix}
{C_{332}}{C_{232}^{-3}} & -{C_{332}}{C_{232}^{-3}} & 0 \\ 
{C_{332}}{C_{232}^{-2}} & 0 & -\frac 12 {C_{332}}{C_{232}^{-2}} \\ 
0 & 0 &  {C^{-1}_{232}}%
\end{pmatrix}%
.
\end{equation*}%
Then $\theta \ast \phi =\big( \triangledown_{13},\triangledown_{33},0\big)
$. So we get the algebra ${\rm M}_{05}.$
\end{itemize}

\item If $\big( C_{131},C_{231},C_{331}\big) =\big( 0,0,1\big),$ then $ C_{232}\neq 0$. Consider the
following matrix:%
\begin{equation*}
\phi=\begin{pmatrix}
0 & {C_{232}^{-2}} & 0 \\ 
1 & -{C_{132}}{C_{232}^{-3}} & -\frac 12 \big(C_{132}+C_{232}C_{332} \big){C_{232}^{-3}}  \\ 
0 & 0 & {C^{-1}_{232}}%
\end{pmatrix}.
\end{equation*}%
Then $\theta \ast \phi = \big( \triangledown_{13},\triangledown_{33},0\big)$. Therefore, there are no new algebras.
\end{itemize}

\noindent
{\bf Second}, we consider algebras with $\dim {\rm A}^{2}=1.$ 
Then let $\big\{e_{1},e_{2},e_{3}\big\}$ be a basis of ${\rm A}$ such that ${\rm A}^{2}=\left\langle e_{1}\right\rangle $. Then the nontrivial multiplications
in ${\rm A}$ is given as follows:%
\begin{longtable}{lcllllll}
${\rm A}$&$
:$&$e_{1}e_{2}=C_{121}e_{1}$ &  $e_{1}e_{3}=C_{131}e_{1}$& $e_{2}e_{2}=C_{221}e_{1}$ 
 & $e_{2}e_{3}=C_{231}e_{1}$ & $e_{3}e_{3}=C_{331}e_{1}$

\end{longtable}\noindent
If $\big(C_{121}, C_{131}\big) = (0,0)$, then $\rm A$ is  associative, 
yielding the algebras ${\rm M}_{01}$ and ${\rm M}_{02}$. 

\medskip

We now consider the case 
$\big(C_{121}, C_{131}\big) \neq (0,0)$. Define 
$\theta\colon {\rm A}\times {\rm A}\rightarrow {\rm A}$ by 

\begin{longtable}{lcl}$\theta$&$
=$&$C_{221}\triangledown_{22}+C_{331}\triangledown
_{33}+C_{121}\triangledown_{12}+C_{131}\triangledown
_{13}+C_{231}\triangledown_{23}$.
\end{longtable} Let $\phi \in {\rm GL}_{3}\big( \mathbb{C}%
\big) $ and write $\theta \ast \phi =\theta ^{\prime }$. Then $\big( 
{\rm A},\theta \big) \cong \big( {\rm A},\theta ^{\prime }\big) 
$. Now, we choose $\phi =\phi _{1}$ where%
\begin{equation*}
\phi _{1}=%
\begin{pmatrix}
1 & 0 & 0 \\ 
0 & a_{22} & a_{23} \\ 
0 & a_{32} & a_{33}%
\end{pmatrix}%
.
\end{equation*}%
Then $\theta ^{\prime }\ =\ C_{221}^{\prime }\triangledown
_{22}+C_{331}^{\prime }\triangledown_{33}+C_{121}^{\prime }\triangledown
_{12}+C_{131}^{\prime }\triangledown_{13}+C_{231}^{\prime }\triangledown
_{23},$ where%
\begin{longtable}{lcl }
$C_{121}^{\prime } $&$=$&$C_{121}a_{22}+C_{131}a_{32},$ \\ 
$C_{131}^{\prime } $&$=$&$C_{121}a_{23}+C_{131}a_{33}.$
\end{longtable}%
Since $\big( C_{121},C_{131}\big) \neq \big( 0,0\big) $, we may assume 
$\big( C_{121}^{\prime },C_{131}^{\prime }\big) =\big( 0,1\big) $.
Next, we choose $\phi =\phi _{2}$ if $C_{221}^{\prime }=0$ or $\phi =\phi
_{3}$ if $C_{221}^{\prime }\neq 0$ where%
\begin{equation*}
\phi _{2}=%
\begin{pmatrix}
1 & -C_{231}^{\prime } & -\frac{1}{2}C_{331}^{\prime } \\ 
0 & 1 & 0 \\ 
0 & 0 & 1%
\end{pmatrix}%
, \ 
\phi _{3}=%
\begin{pmatrix}
1 & -C_{231}^{\prime}( C_{221}^{\prime })^{-\frac{1}{2}} & -\frac{1}{2}%
C_{331}^{\prime } \\ 
0 &  ( {C_{221}^{\prime }})^{-\frac{1}{2}} & 0 \\ 
0 & 0 & 1%
\end{pmatrix}%
.
\end{equation*}%
Then $\theta ^{\prime }\ast \phi \in \big\{ \triangledown
_{13},\triangledown_{13}+\triangledown_{22}\big\} $. Then we get the
algebras ${\rm M}_{04}^{0}$ and ${\rm M}_{07}$.
\end{proof}

\subsection{Derived commutative associative  algebras}\label{derasscom}

\begin{definition}
A commutative algebra ${\rm A}$ is called a derived commutative associative algebra if the derived algebra $
{\rm A}^{2}$ is a  commutative associative algebra.
\end{definition}

It is clear that any metabelian  commutative  algebra is a derived commutative associative 
 algebra. The variety of derived  commutative associative algebras contains all  commutative associative algebras.
Let ${\rm A}$ be a derived  commutative associative algebra. 
If either ${\rm A}^{2}=%
{\rm A}$ or ${\rm A}^{2}=0$, then ${\rm A}$\ is a  commutative associative  algebra.

\begin{theoremA2}\label{thm:alg_derComAs}
Let ${\rm A}$ be a complex $3$-dimensional derived commutative associative algebra. Then ${\rm A}$
is 
a metabelian commutative algebra listed in {\rm Theorem A1}, 
or a commutative associative algebra listed in {\rm Proposition~\ref{3-dim assoc}}, 
or it is
isomorphic to one of the following algebras:

\begin{longtable}{lllllllllll}
${\rm A}_{01}^{\alpha, \beta, \gamma}$ & $:$ & $e_{1}e_{1} = e_{1}$ & $e_{2}e_{2} = e_{2}$ & $e_{1}e_{3} = \alpha e_{2}$ \\&& $e_{2}e_{3} = \beta e_{1}$ & $e_{3}e_{3} = e_{1} + \gamma e_{2}$ \\ 
${\rm A}_{02}^{\alpha}$ & $:$ & $e_{1}e_{1} = e_{1}$ & $e_{2}e_{2} = e_{2}$ & $e_{1}e_{3} = e_{2}$ & $e_{2}e_{3} = \alpha e_{1}$ & \\ 
${\rm A}_{03}^{\alpha}$ & $:$ & $e_{1}e_{1} = e_{1}$ & $e_{2}e_{3} = e_{2}$ & $e_{3}e_{3} = \alpha e_{1}$ & & \\ 
${\rm A}_{04}^{\alpha}$ & $:$ & $e_{1}e_{1} = e_{1}$ & $e_{2}e_{3} = e_{1} + e_{2}$ & $e_{3}e_{3} = \alpha e_{1}$ & & \\ 
${\rm A}_{05}^{\alpha, \beta}$ & $:$ & $e_{1}e_{1} = e_{1}$ & $e_{1}e_{3} = e_{2}$ & $e_{2}e_{3} = \alpha e_{1} + e_{2}$ & $e_{3}e_{3} = \beta e_{1}$ & \\ 
${\rm A}_{06}^{\alpha}$ & $:$ & $e_{1}e_{1} = e_{1}$ & $e_{1}e_{3} = e_{2}$ & $e_{2}e_{3} = e_{1}$ & $e_{3}e_{3} = \alpha e_{2}$ & \\ 
${\rm A}_{07}$ & $:$ & $e_{1}e_{1} = e_{1}$ & $e_{1}e_{3} = e_{2}$ & & & \\ 
${\rm A}_{08}$ & $:$ & $e_{1}e_{1} = e_{1}$ & $e_{1}e_{3} = e_{2}$ & $e_{3}e_{3} = e_{2}$ & & \\ 
${\rm A}_{09}^{\alpha}$ & $:$ & $e_{1}e_{1} = e_{1}$ & $e_{1}e_{3} = e_{2}$ & $e_{3}e_{3} = e_{1} + \alpha e_{2}$ & & \\ 
${\rm A}_{10}$ & $:$ & $e_{1}e_{1} = e_{1}$ & $e_{2}e_{3} = e_{1}$ & $e_{3}e_{3} = e_{2}$ & & \\ 
${\rm A}_{11}$ & $:$ & $e_{1}e_{1} = e_{1}$ & $e_{3}e_{3} = e_{1} + e_{2}$ & & & \\ 
${\rm A}_{12}^{\alpha, \beta}$ & $:$ & $e_{1}e_{1} = e_{1}$ & $e_{1}e_{2} = e_{2}$ & $e_{1}e_{3} = e_{1}$  \\&& $e_{2}e_{3} = \alpha e_{1}$ & $e_{3}e_{3} = \beta e_{1} + e_{2}$ \\ 
${\rm A}_{13}^{\alpha}$ & $:$ & $e_{1}e_{1} = e_{1}$ & $e_{1}e_{2} = e_{2}$ & $e_{1}e_{3} = e_{1}$  \\&& $e_{2}e_{3} = e_{1}$ & $e_{3}e_{3} = \alpha e_{1}$ \\ 
${\rm A}_{14}^{\alpha}$ & $:$ & $e_{1}e_{1} = e_{1}$ & $e_{1}e_{2} = e_{2}$ & $e_{1}e_{3} = e_{1}$ & $e_{3}e_{3} = \alpha e_{1}$ & \\ 
${\rm A}_{15}^{\alpha}$ & $:$ & $e_{1}e_{1} = e_{1}$ & $e_{1}e_{2} = e_{2}$ & $e_{2}e_{3} = e_{1}$ & $e_{3}e_{3} = e_{1} + \alpha e_{2}$ & \\ 
${\rm A}_{16}$ & $:$ & $e_{1}e_{1} = e_{1}$ & $e_{1}e_{2} = e_{2}$ & $e_{3}e_{3} = e_{1} + e_{2}$ & & \\ 
${\rm A}_{17}$ & $:$ & $e_{1}e_{1} = e_{1}$ & $e_{1}e_{2} = e_{2}$ & $e_{3}e_{3} = e_{1}$ & & \\ 
${\rm A}_{18}$ & $:$ & $e_{1}e_{1} = e_{1}$ & $e_{1}e_{2} = e_{2}$ & $e_{2}e_{3} = e_{1}$ & $e_{3}e_{3} = e_{2}$ & \\ 
${\rm A}_{19}$ & $:$ & $e_{1}e_{1} = e_{1}$ & $e_{1}e_{2} = e_{2}$ & $e_{2}e_{3} = e_{1}$ & & \\ 
${\rm A}_{20}$ & $:$ & $e_{1}e_{1} = e_{1}$ & $e_{1}e_{2} = e_{2}$ & $e_{3}e_{3} = e_{2}$ & & \\ 
${\rm A}_{21}^{\alpha}$ & $:$ & $e_{1}e_{1} = e_{2}$ & $e_{2}e_{3} = e_{1} + e_{2}$ & $e_{3}e_{3} = \alpha e_{2}$ & & \\ 
${\rm A}_{22}$ & $:$ & $e_{1}e_{1} = e_{2}$ & $e_{2}e_{3} = e_{1}$ & & & \\ 
${\rm A}_{23}$ & $:$ & $e_{1}e_{1} = e_{2}$ & $e_{2}e_{3} = e_{1}$ & $e_{3}e_{3} = e_{2}$ & & \\ 
${\rm A}_{24}^{\alpha}$ & $:$ & $e_{1}e_{1} = e_{2}$ & $e_{1}e_{3} = e_{1}$ & $e_{2}e_{3} = \alpha e_{2}$ & & \\ 
${\rm A}_{25}$ & $:$ & $e_{1}e_{1} = e_{2}$ & $e_{1}e_{3} = e_{1}$ & $e_{3}e_{3} = e_{2}$ & & \\ 
${\rm A}_{26}$ & $:$ & $e_{1}e_{1} = e_{2}$ & $e_{1}e_{3} = e_{1} + e_{2}$ & $e_{2}e_{3} = e_{2}$ & & \\ 
${\rm A}_{27}$ & $:$ & $e_{1}e_{1} = e_{2}$ & $e_{3}e_{3} = e_{1}$ & & & \\ 
${\rm A}_{28}$ & $:$ & $e_{1}e_{1} = e_{2}$ & $e_{2}e_{3} = e_{2}$ & $e_{3}e_{3} = e_{1}$ & & \\ 
${\rm A}_{29}$ & $:$&$e_{1}e_{1}=e_{1}$&$e_{2}e_{3}=e_{1}$\\
${\rm A}_{30}$ & $:$&$e_{1}e_{1}=e_{1}$&$e_{3}e_{3}=e_{1}$
\end{longtable}
\noindent
All listed algebras are non-isomorphic, excepting:
\begin{center}
$
{\rm A}_{01}^{\alpha ,\beta ,\gamma }\cong 
{\rm A}_{01}^{-\alpha,\beta,\gamma}\cong 
{\rm A}_{01}^{\alpha,-\beta,\gamma}\cong 
{\rm A}_{01}^{\beta \gamma^{-\frac 12}, \alpha \gamma^{-\frac 12}, \gamma^{-1}}, \
{\rm A}_{02}^{\alpha }\cong {\rm A}_{02}^{\alpha^{-1}}, \ 
{\rm A}_{06}^{\alpha }\cong {\rm A}_{06}^{-\alpha},\ 
{\rm A}_{09}^{\alpha }\cong {\rm A}_{09}^{-\alpha},\ 
{\rm A}_{15}^{\alpha }\cong {\rm A}_{15}^{-\alpha}.$
\end{center}

\end{theoremA2}

\begin{proof}
{\bf Firstly}, we consider the case when $\dim {\rm A}^{2} = 2.$
We may assume 
\[
{\rm A}^{2} = \langle e_{1}, e_{2} \rangle \in 
\big\{ {\mathfrak J}_{01}, {\mathfrak J}_{02}, {\mathfrak J}_{03}, {\mathfrak J}_{05} \big\},
\]
where ${\mathfrak J}_{04}$ is excluded, as it is non-associative. Extend a basis $\left\{ e_{1},e_{2}\right\} $ of ${\rm A}^{2}$
to a basis $\left\{ e_{1},e_{2},e_{3}\right\} $ of ${\rm A}$. Since our aim is to classify non-associative algebras ${\rm A}$, we only consider those cases where the multiplication on ${\rm A}$ is not associative. The associative algebras among these possibilities have already been described in Proposition~\ref{3-dim assoc}. Hence, in what follows, we assume that ${\rm A}$ is non-associative.

\begin{enumerate}[I.]
\item
\underline{${\rm A}^{2}={\mathfrak J}_{01}$.} Then the multiplication
table of ${\rm A}$ is defined as follows: 
\begin{longtable}{lcllllll}
${\rm A}$&$:$& $e_{1}e_{1}=e_{1} $&
$e_{1}e_{3}=C_{131}e_{1}+C_{132}e_{2}$&
$e_{2}e_{3}=C_{231}e_{1}+C_{232}e_{2}$\\
&&& $ e_{2}e_{2}=e_{2} $&
$e_{3}e_{3}=C_{331}e_{1}+C_{332}e_{2}$
\end{longtable}
Without loss of generality, we may assume $C_{131} = C_{232} = 0$ via the change of basis:
\[
e_1' = e_1, \quad e_2' = e_2, \quad e_3' = -C_{131} e_1 - C_{232} e_2 + e_3.
\]

Further if $C_{131} = C_{232} = 0$, then $\mathrm{A}$ is a commutative associative algebra precisely when
\[
C_{331} = C_{332} = C_{132} = C_{231} = 0.
\]
Hence we assume $\big(C_{331}, C_{332}, C_{132}, C_{231}\big) \neq (0, 0, 0, 0)$. Now, let us consider the following cases:

\begin{itemize}
\item If $C_{331}\neq 0,$ then with the change of basis: 
$e_{1}^{\prime}=e_{1},$ \ 
$e_{2}^{\prime }=e_{2},$ \ 
$e_{3}^{\prime }={C^{-\frac{1}{2}}_{331}}e_{3},$
we get the family of algebras ${\rm A}_{01}^{\alpha ,\beta ,\gamma }$. Moreover,
\begin{center}$%
{\rm A}_{01}^{\alpha ,\beta ,\gamma }\cong {\rm A}_{01}^{\alpha
^{\prime},\beta^{\prime},\gamma ^{\prime }}$ if 
$\alpha ^{2}=\alpha^{\prime 2},$ \ 
$\beta ^{2}=\beta ^{\prime 2},$ \ 
$\gamma =\gamma ^{\prime }$ or $\alpha ^{2}=\gamma \beta ^{\prime 2},$\ 
$\beta ^{2}=\gamma \alpha ^{\prime2},$ \ $\gamma \gamma ^{\prime }=1$.\end{center}

\item If $C_{331}=0$ and $C_{332}\neq 0,$ then with the change of basis: 
\begin{center}
$e_{1}^{\prime }=e_{2},$ \  
$e_{2}^{\prime }=e_{1},$ \ 
$e_{3}^{\prime }=C^{-\frac{1}{2}}_{332}e_{3},$ 
\end{center}
we get the family of algebras ${\rm A}_{01}^{\alpha ,\beta ,0}$.

\item If $C_{331}=C_{332}=0$ and $C_{132}\neq 0,$ then after changing 
\begin{center}
    $e_{1}^{\prime}=e_{1},$ \ 
$e_{2}^{\prime }=e_{2},$ \ 
$e_{3}^{\prime }= {C^{-1}_{132}}e_{3},$ 
\end{center} we get
the family of algebras ${\rm A}_{02}^{\alpha }$. Furthermore, ${\rm A}_{02}^{\alpha }\cong {\rm A}_{02}^{\beta}$ if 
$ \alpha =  \beta ^{-1}$.

\item If $C_{331}=C_{332}= C_{132}=0,$ then after changing 
$e_{1}^{\prime}=e_{2},$ \ 
$e_{2}^{\prime }=e_{1},$ \ 
$e_{3}^{\prime }= {C^{-1}_{231}}e_{3},$ we get
the algebra ${\rm A}_{02}^{0}$.
\end{itemize}

\item
\underline{${\rm A}^{2}={\mathfrak J}_{02}$.} The nontrivial
multiplications in ${\rm A}$ is given as follows: 
\begin{longtable}{lcllll}
${\rm A}$ & $:$&
$e_{1}e_{1}=e_{1} $&$e_{1}e_{3}=C_{131}e_{1}+C_{132}e_{2} $&$e_{2}e_{3}=C_{231}e_{1}+C_{232}e_{2}$\\
&&&&$e_{3}e_{3}=C_{331}e_{1}+C_{332}e_{2}$
\end{longtable}%
Applying the change of basis:
$e_{1}'=e_{1},$ 
$e_{2}'=e_{2},$ 
$e_{3}'=-C_{131}e_{1}+e_{3},$
we may assume $C_{131}=0$. In this reduced form, the algebra ${\rm A}$ is commutative and associative precisely when
\[
C_{331}=C_{132}=C_{231}=C_{232}=0.
\]
Since $\dim {\rm A}^{2}=2$, it follows that
$\big(C_{232},\, C_{132},\, C_{332} \big)\neq (0,0,0).$

\begin{itemize}
\item If $C_{232}\neq 0,$ then, with the change of basis: 
\begin{center}
    $e_{1}^{\prime}=e_{1},$ \ 
$e_{2}^{\prime }=e_{2},$ \ 
$e_{3}^{\prime }=-\frac{C_{332}}{2C_{232}^{2}}%
e_{2}+\frac{1}{C_{232}}e_{3},$ 
\end{center} we can make $C_{232}=1$ and $C_{332}=0$. So we
assume $C_{232}=1$ and $C_{332}=0$.

\begin{itemize}
\item If $C_{132}=C_{231}=0$, we get the family of algebras ${\rm A}_{03}^{\alpha }.$

\item If $C_{132}=0$ and $C_{231}\neq 0,$ 
then after changing  
$e_{1}^{\prime }=e_{1},$ \ 
$e_{2}^{\prime }= C^{-1}_{231}e_{2},$ \ 
$e_{3}^{\prime}=e_{3},$ we get the family of algebras ${\rm A}_{04}^{\alpha }$. 

\item If $C_{132}\neq 0,$
then with the change of basis: 
$e_{1}^{\prime}=e_{1},$ \ 
$e_{2}^{\prime }=C_{132}e_{2},$ \ 
$e_{3}^{\prime }=e_{3},$ we get the family of
algebras ${\rm A}_{05}^{\alpha ,\beta }$. 
\end{itemize}

\item If $C_{232}=0$ and $C_{132}\neq 0,$ 
then with the change of basis: 
$e_{1}^{\prime}=e_{1},$ \ 
$e_{2}^{\prime }=C_{132}e_{2},$ \ 
$e_{3}^{\prime }=e_{3},$ 
we can make $%
C_{132}=1$. So we assume $C_{132}=1$.

\begin{itemize}
\item If $C_{231}\neq 0$, then with the change of basis: 
\begin{center}
    $e_{1}^{\prime}=e_{1},$ \ 
$e_{2}^{\prime }={C^{-\frac{1}{2}}_{231}}e_{2},$ \ 
$e_{3}^{\prime }=-\frac{1}{2} {C_{331}}{C_{231}^{-\frac{3}{2}}}e_{2}+{{C^{-\frac{1}{2}}_{231}}}e_{3},$
\end{center}
 we get the family of algebras ${\rm A}_{06}^{\alpha }$.
Moreover,  ${\rm A}_{06}^{\alpha }\cong {\rm A}_{06}^{\beta }$ if and
only if $\alpha =- \beta$.

\item If $C_{231}=C_{331}=C_{332}=0$, we get the algebra ${\rm A}%
_{07}$. 

\item If $C_{231}=C_{331}=0$ and $C_{332}\neq 0$, 
then by  changing \begin{center}
$e_{1}^{\prime }=e_{1},$ \ 
$e_{2}^{\prime }={C^{-1}_{332}}e_{2},$ \ 
$e_{3}^{\prime}={C^{-1}_{332}}e_{3},$ 
\end{center}we get the algebra ${\rm A}_{08}$.

\item If $C_{231}=0$ and $C_{331}\neq 0,$ 
then  by changing
\begin{center}
$e_{1}^{\prime }=e_{1},$ \ 
$e_{2}^{\prime }= {C^{-\frac{1}{2}}_{331}}e_{2},$ \ 
$e_{3}^{\prime }={C^{-\frac{1}{2}}_{331}}e_{3},$
\end{center} we get the family of algebras ${\rm A}_{09}^{\alpha }$. Moreover, ${\rm A}_{09}^{\alpha }\cong {\rm A}_{09}^{\beta }$ if and only if $\alpha =-\beta$.

\end{itemize}

\item If $C_{232}=C_{132}=0,$ then $C_{332}\neq 0$ since $\dim {\rm A}^{2}=2$. Moreover,
we have $\big( C_{231},C_{331}\big) \neq \left( 0,0\right) $ since
otherwise ${\rm A}$ is a commutative associative algebra. Then with the change of basis: 
$e_{1}^{\prime }=e_{1},$ \ 
$e_{2}^{\prime }=e_{2},$ \ 
$e_{3}^{\prime }={C^{-\frac{1}{2}}_{332}}e_{3}$ we can make $C_{332}=1$. So we assume $C_{332}=1$.
\begin{itemize}
    \item If $%
C_{231}\neq 0$, then with the change of basis: 
\begin{center}
    $e_{1}^{\prime}=e_{1},$ \ 
    $e_{2}^{\prime }={C^{-\frac{2}{3}}_{231}}e_{2},$ \ 
    $e_{3}^{\prime}=-\frac{1}{2}C_{231}^{-\frac{4}{3}}C_{331}e_{2}+{C^{-\frac{1}{3}}_{231}}e_{3},$
\end{center} we get the algebra ${\rm A}_{10}$. 

  \item If $C_{231}=0$, then $%
C_{331}\neq 0$ and with the change of basis: 
\begin{center}
    $e_{1}^{\prime}=e_{1},$ \ 
${e_{2}^{\prime }={C^{-1}_{331}}e_{2}},$ \ 
$e_{3}^{\prime }={C^{-\frac{1}{2}}_{331}}e_{3},$ 
\end{center}we get the algebra ${\rm A}_{11}$.

\end{itemize}
\end{itemize}

\item
\underline{${\rm A}^{2}={\mathfrak J}_{03}$.} Then the nontrivial
multiplications in ${\rm A}$ is given as follows: 
\begin{longtable}{lclllll}
${\rm A}$ & $:$ & $e_{1}e_{1}=e_{1}$ & $e_{1}e_{3}=C_{131}e_{1}+C_{132}e_{2}$ & $e_{2}e_{3}=C_{231}e_{1}+C_{232}e_{2}$ \\
&&& $e_{1}e_{2}=e_{2}$ &$ e_{3}e_{3}=C_{331}e_{1}+C_{332}e_{2}$
\end{longtable}%
Then with the change of basis: 
$e_{1}^{\prime }=e_{1},$ \ 
$e_{2}^{\prime}=e_{2},$   \ 
$e_{3}^{\prime }=-C_{232}e_{1}-C_{132}e_{2}+e_{3},$ we can make $
C_{132}=C_{232}=0$. So we assume $C_{132}=C_{232}=0$. Hence the nontrivial
multiplications in ${\rm A}$ is given as follows:%
\begin{longtable}{lclllll}
${\rm A}$ &
$:$ & $e_{1}e_{1}=e_{1}$ & $e_{1}e_{2}=e_{2}$ & $e_{1}e_{3}=C_{131}e_{1}$ & $e_{2}e_{3}=C_{231}e_{1} $& $e_{3}e_{3}=C_{331}e_{1}+C_{332}e_{2}$
\end{longtable}%
So ${\rm A}$ is a commutative associative algebra if $C_{131}=C_{231}=C_{331}=C_{332}=0$.

\begin{itemize}
\item We consider the case $C_{131}\neq 0$.
\begin{itemize}
\item If $C_{332}\neq 0$, then by  changing  
$e_{1}^{\prime}=e_{1},$ \ 
$e_{2}^{\prime }={C_{332}}{C_{131}^{-2}}e_{2},$ \ 
$e_{3}^{\prime }={C^{-1}_{131}}e_{3},$ we get the family of algebras ${\rm A}_{12}^{\alpha ,\beta
}$. 

\item If $C_{332}=0$ and $C_{231}\neq 0,$ then by  changing  
\begin{center}$e_{1}^{\prime }=e_{1},$ \ 
$e_{2}^{\prime }={C_{131}}{C^{-1}_{231}}e_{2},$ \ 
${e_{3}^{\prime }={C^{-1}_{131}}e_{3}},$ 
\end{center} we get the family of algebras $\rm{A}_{13}^{\alpha }$. 

\item If $C_{332}= C_{231}=0,$ then by  changing   
$e_{1}^{\prime}=e_{1},$ \ 
$e_{2}^{\prime }=e_{2},$ \ 
$e_{3}^{\prime }= {C^{-1}_{131}}e_{3},$ we get
the family of algebras ${\rm A}_{14}^{\alpha }$. 
\end{itemize}

\item Now, we consider the second case when $C_{131}=0$ and $C_{331}\neq 0$.

\begin{itemize}
\item If $C_{231}\neq 0,$ then by  changing  
$e_{1}^{\prime}=e_{1},$ \ 
$e_{2}^{\prime }={C^{\frac{1}{2}}_{331}}{C^{-1}_{231}}e_{2},$ \ 
$e_{3}^{\prime }= {C^{-\frac{1}{2}}_{331}}e_{3},$ we get the family of algebras 
${\rm A}_{15}^{\alpha}$. Moreover, ${\rm A}_{15}^{\alpha }\cong {\rm A}_{15}^{\beta }$ if
and only if $\alpha =-\beta$.

\item If $C_{231}=0$ and $C_{332}\neq 0,$ 
then by  changing  
\begin{center}
$e_{1}^{\prime}=e_{1},$ \ 
$e_{2}^{\prime }={C_{332}}{C^{-1}_{331}}e_{2},$ \ 
${e_{3}^{\prime }= {C^{-\frac{1}{2}}_{331}}e_{3}},$ 
\end{center} we get the algebra ${\rm A}_{16}$.

\item If $C_{231}=C_{332}=0,$ then by  changing   
    $e_{1}^{\prime}=e_{1},$ \ 
$e_{2}^{\prime }=e_{2},$ \ 
$e_{3}^{\prime }={C^{-\frac{1}{2}}_{331}}e_{3},$
we get the algebra ${\rm A}_{17}$.
\end{itemize}

\item Finally, we consider the third case when $C_{131}=0$ and $C_{331}=0$.

\begin{itemize}
\item If $C_{231}C_{332}\neq 0,$ then with the change of basis: 
\begin{center}
$e_{1}^{\prime}=e_{1},$\ 
$e_{2}^{\prime }= C_{332}^{\frac 13}C_{231}^{-\frac 23}e_{2},$ \ 
$e_{3}^{\prime }= C^{-\frac 13}_{231}C^{-\frac 13}_{332}e_{3},$
\end{center} we get the
algebra ${\rm A}_{18}$.

\item If $C_{231}\neq 0$ and $C_{332}=0,$ then with the change of basis: 
\begin{center}
$e_{1}^{\prime }=e_{1},$ \ 
$e_{2}^{\prime }=e_{2},$ \ 
$e_{3}^{\prime }={C^{-1}_{231}}e_{3},$
\end{center} we get the algebra ${\rm A}_{19}$.

\item If $C_{231}=0,$ then since ${\rm A}$ is non-associative, we have $C_{332}\neq 0$. Then, with the change of basis: 
$e_{1}^{\prime }=e_{1},$ \ 
$e_{2}^{\prime}=C_{332}e_{2},$ \ 
$e_{3}^{\prime }=e_{3},$ we get the algebra ${\rm A}_{20}$.
\end{itemize}
\end{itemize}

\item
\underline{${\rm A}^{2}={\mathfrak J}_{05}$.} Then the nontrivial
multiplications in ${\rm A}$ is given as follows: 
\begin{longtable}{lclllll}
${\rm A}$&$:$&
$e_{1}e_{1}=e_{2}$ & $e_{1}e_{3}=C_{131}e_{1}+C_{132}e_{2}$ & $e_{2}e_{3}=C_{231}e_{1}+C_{232}e_{2}$ \\&&&& $e_{3}e_{3}=C_{331}e_{1}+C_{332}e_{2}$
\end{longtable}
Then with the change of basis: 
$e_{1}^{\prime }=e_{1},$ \ 
$e_{2}^{\prime}=e_{2},$ \ 
$e_{3}^{\prime }=-C_{132}e_{1}+e_{3},$ 
we can make $C_{132}=0$. So we
assume $C_{132}=0$. Moreover, ${\rm A}$ is associative if 
\begin{center}$
C_{131}=C_{231}=C_{331}=C_{232}=0$.\end{center} Since $\dim {\rm A}^{2}=2$, we have  $\big(
C_{131},C_{231},C_{331}\big) \neq \big( 0,0,0\big) $.

\begin{itemize}
\item Case $C_{231}\neq 0$. If we consider the change of basis: 
\begin{center}

$e_{1}^{\prime}=e_{1}-\frac{C_{131}}{C_{231}}e_{2},$ \ 
$e_{2}^{\prime }=e_{2},$\ 
$e_{3}^{\prime }=\frac{C_{131}C_{232}}{C_{231}^{2}}e_{1}-\frac{ 
2C_{232}C_{131}^{2}+C_{231}C_{331}  }{2C_{231}^{3}}e_{2}+\frac{1}{%
C_{231}}e_{3},$

\end{center}we can make $C_{231}=1$ and $C_{131}=C_{331}=C_{132}=0$.
Whence, the nontrivial multiplications in ${\rm A}$ is given as follows: 
\begin{longtable}{lcllll}
${\rm A}$&$:$&$e_{1}e_{1}=e_{2}$&$e_{2}e_{3}=e_{1}+C_{232}e_{2}$&$e_{3}e_{3}=C_{332}e_{2}$
\end{longtable}

\begin{itemize}
\item If $C_{232}\neq 0,$ then by changing  
$e_{1}^{\prime}=C_{232}e_{1},$ \ 
$e_{2}^{\prime }=C_{232}^{2}e_{2},$ \ 
${e_{3}^{\prime }={C^{-1}_{232}}e_{3}},$ we get the family of algebras ${\rm A}_{21}^{\alpha }$. 

\item If $C_{232}=C_{332}=0$, then we get the algebra ${\rm A}_{22}$. 

\item If $C_{232}=0$ and $C_{332}\neq 0$, then with the change of basis: 
\begin{center} 
$e_{1}^{\prime}={C^{\frac{1}{4}}_{332}}e_{1},$ \ 
${e_{2}^{\prime}={C^{\frac{1}{2}}_{332}}e_{2}},$ \ 
${e_{3}^{\prime}={C^{-\frac{1}{4}}_{332}}e_{3}},$
\end{center} we get the algebra ${\rm A}_{23}$.
\end{itemize}

\item Case $C_{231}=0$ and
 $C_{131}\neq 0$. Then, with the change of basis: 
\begin{center}
    $e_{1}^{\prime}=e_{1},$ \ 
$e_{2}^{\prime }=e_{2},$ \ 
${e_{3}^{\prime }=-\frac{1}{2} C_{331} C_{131}^{-2}e_{1}+ {C^{-1}_{131}}e_{3}},$ 
\end{center}
we may assume $C_{131}=1$ and $C_{331}=0$. Whence,
the nontrivial multiplications in ${\rm A}$ is given as follows: 
\begin{longtable}{lclllll}
${\rm A}$&$
:$&$e_{1}e_{1}=e_{2}$&$e_{1}e_{3}=e_{1}+C_{132}e_{2}$&$e_{2}e_{3}=C_{232}e_{2}$&$e_{3}e_{3}=C_{332}e_{2}$ 
\end{longtable}
\begin{itemize}
\item If $C_{232}\notin \big\{ 0,1\big\},$ then with the change of basis: 
\begin{center}
$e_{1}^{\prime}=e_{1}-\frac{C_{132}}{C_{232}-1}e_{2},$\ 
$e_{2}^{\prime }=e_{2},$\ 
${e_{3}^{\prime}=-\frac{1}{2} {C_{332}}{C^{-1}_{232}}e_{2}+e_{3}},$
\end{center} we get the algebras ${\rm A}
_{24}^{\alpha \notin\{ 0,1\}}$. Moreover, ${\rm A}_{24}^{\alpha }\cong 
{\rm A}_{24}^{\beta }$ if and only if $\alpha =\beta $.

\item If $C_{232}=C_{332}=0$, then with the change of basis: 
\begin{center}$e_{1}^{\prime }=e_{1}+C_{132}e_{2},$\ 
$e_{2}^{\prime }=e_{2},$ \ 
$e_{3}^{\prime}=e_{3},$ 
\end{center} we get the algebra ${\rm A}_{24}^{0}$. 

\item If $C_{232}=0$ and $C_{332}\neq0 $, then by changing
\begin{center}${e_{1}^{\prime }={C^{\frac{1}{2}}_{332}}e_{1}+C_{132}{C^{\frac{1}{2}}_{332}}e_{2}},$ \ 
${e_{2}^{\prime }=C_{332}e_{2}},$ \ 
${e_{3}^{\prime}=e_{3}},$
\end{center} we get the algebra ${\rm A}_{25}$.

\item If $C_{232}=1$ and $C_{132}=0$,
then with the change of basis: 
\begin{center}
${e_{1}^{\prime }=e_{1}},$ \ 
${e_{2}^{\prime }=e_{2}},$ \ 
${e_{3}^{\prime }=-\frac{1}{2}C_{332}e_{2}+e_{3}},$ 
\end{center}we get the algebra ${\rm A}_{24}^{ 1}$.

\item If $C_{232}=1$ and  $%
C_{132}\neq 0$,
then with the change of basis: 
\begin{center}$e_{1}^{\prime}=C_{132}e_{1},$ \ 
${e_{2}^{\prime }=C_{132}^{2}e_{2}},$ \ 
${e_{3}^{\prime }=-\frac{1}{2}C_{332}e_{2}+e_{3}},$ 
\end{center}
we get the algebra ${\rm A}_{26}$.
\end{itemize}

\item If $C_{231}=C_{131}=0,$ then $C_{331}\neq 0$ and we consider the following subcases.
\begin{itemize}
    \item If $C_{232}=0$, then with the
change of basis: \begin{center}
${e_{1}^{\prime }=C_{331}e_{1}+C_{332} e_{2}},$ \ 
${e_{2}^{\prime}=C_{331}^{2}e_{2}},$ \ 
${e_{3}^{\prime }=e_{3}},$ 
\end{center}we get the algebra
${\rm A}_{27}$. 

\item If $C_{232}\neq 0$, then with the change of basis: 
\begin{center}
$e_{1}^{\prime }={C_{331}}{C_{232}^{-2}}e_{1},$\ 
$e_{2}^{\prime }={C_{331}^{2}}{C_{232}^{-4}}e_{2},$ \ 
$e_{3}^{\prime }=-\frac{1}{2} {C_{332}}{C_{232}^{-2}}e_{2}+C^{-1}_{232}e_{3},$
\end{center}we get the algebra ${\rm A}_{28}$.
\end{itemize}
 
\end{itemize}

\end{enumerate}

\medskip 
\noindent
{\bf Secondly}, we consider the case when $\dim {\rm A}^{2} = 1$
 and ${\rm A}^{2}{\rm A}^{2}\neq 0$ $\big($i.e., $\rm A$ is not metabelian$\big).$
We may assume ${\rm A}^{2}=\left\langle e_{1}\right\rangle $ such that $%
e_{1}e_{1}=e_{1}$. Extend a basis $\left\{ e_{1}\right\} $ of ${\rm A}
^{2}$ to a basis $\left\{ e_{1},e_{2},e_{3}\right\} $ of ${\rm A}$. Then
the nontrivial multiplications in ${\rm A}$ is given as follows:%
\begin{longtable}{lcllll}
${\rm A}$&$:$&$e_{1}e_{1}=e_{1}$&$e_{2}e_{2}=C_{221}e_{1}$&
$e_{3}e_{3}=C_{331}e_{1}$\\&&$e_{1}e_{2}=C_{121}e_{1}$ &$e_{1}e_{3}=C_{131}e_{1}$&$e_{2}e_{3}=C_{231}e_{1}$
\end{longtable}
Then ${\rm A}$ is associative if the following relations hold:%
\begin{longtable}{lcllcllcl}
$C_{221}$&$=$&$C_{121}^{2},$&
$C_{331}$&$=$&$C_{131}^{2},$& 
$C_{121}C_{131}$&$=$&$C_{231}.$
\end{longtable}%
Define $\theta\colon \rm{A}\times {\rm A}\rightarrow {\rm A}$ by
\begin{longtable}{lcl}
$\theta$&$=$&$\triangledown_{11}+C_{221}\triangledown_{22}+C_{331}\triangledown
_{33}+C_{121}\triangledown_{12}+C_{131}\triangledown
_{13}+C_{231}\triangledown_{23}$.
\end{longtable} 
Let $\phi \in {\rm GL}_{3}\big( \mathbb{C}%
\big) $ and write $\theta \ast \phi =\theta ^{\prime }$. Then $\big( 
{\rm A},\theta \big) \cong \big( {\rm A},\theta ^{\prime }\big) 
$. Now, we choose $\phi =\phi _{1}$ where:%
\begin{equation*}
\phi _{1}=%
\begin{pmatrix}
1 & -C_{121} & -C_{131} \\ 
0 & 1 & 0 \\ 
0 & 0 & 1%
\end{pmatrix}%
,
\end{equation*}%
then $\theta ^{\prime }=\triangledown_{11}+C_{221}^{\prime }\triangledown
_{22}+C_{331}^{\prime }\triangledown_{33}+C_{231}^{\prime }\triangledown
_{23},$ where 
\begin{center}
$C_{221}^{\prime }=C_{221}-C_{121}^{2},$ \ 
$C_{331}^{\prime}=C_{331}-C_{131}^{2},$ \  
$C_{231}^{\prime }=C_{231}-C_{121}C_{131}$.\end{center} Since $%
{\rm A}$ is non-associative, 
$\big( C_{221}^{\prime
},C_{331}^{\prime },C_{231}^{\prime }\big) \neq \big( 0,0,0\big) $.
Next we choose $\phi =\phi _{2}$ where:%
\begin{equation*}
\phi _{2}=%
\begin{pmatrix}
1 & 0 & 0 \\ 
0 & a_{22} & a_{23} \\ 
0 & a_{32} & a_{33}%
\end{pmatrix}%
.
\end{equation*}%
Then $\theta ^{\prime }\ast \phi =\triangledown_{11}+C_{221}^{\prime
\prime }\triangledown_{22}+C_{331}^{\prime \prime }\triangledown
_{33}+C_{231}^{\prime \prime }\triangledown_{23},$ where%
\begin{longtable}{lcl}
$C_{221}^{\prime \prime } $&$=$&$C_{221}^{\prime }a_{22}^{2}+2C_{231}^{\prime
}a_{22}a_{32}+C_{331}^{\prime }a_{32}^{2},$ \\
$C_{331}^{\prime \prime } $&$=$&$C_{221}^{\prime }a_{23}^{2}+2C_{231}^{\prime
}a_{23}a_{33}+C_{331}^{\prime }a_{33}^{2},$ \\
$C_{231}^{\prime \prime } $&$=$&$C_{221}^{\prime }a_{22}a_{23}+C_{231}^{\prime
}a_{22}a_{33}+C_{231}^{\prime }a_{23}a_{32}+C_{331}^{\prime }a_{32}a_{33}.$
\end{longtable}%
\noindent Whence, 
\begin{longtable}{lcl}$
\begin{pmatrix}
C_{231}^{\prime \prime } & C_{331}^{\prime \prime } \\ 
-C_{221}^{\prime \prime } & -C_{231}^{\prime \prime }%
\end{pmatrix} $&$
= $&$
\begin{pmatrix}
a_{33} & -a_{23} \\ 
-a_{32} & a_{22}%
\end{pmatrix}%
\begin{pmatrix}
C_{231}^{\prime } & C_{331}^{\prime } \\ 
-C_{221}^{\prime } & -C_{231}^{\prime }%
\end{pmatrix}%
\begin{pmatrix}
a_{22} & a_{23} \\ 
a_{32} & a_{33}%
\end{pmatrix}%
$.
\end{longtable}\noindent
So, we may assume 
$\big( C_{221}^{\prime \prime },C_{331}^{\prime \prime
},C_{231}^{\prime \prime }\big) \in \big\{ \big( 0,0,1\big) ,\big(
0,1,0\big) \big\}$, which leads to exactly two non-isomorphic algebras, denoted by ${\rm A}_{29}$ and ${\rm A}
_{30}$.

\end{proof}

\subsection{Derived Jordan algebras}\label{derJ}

\begin{definition}
A commutative algebra ${\rm A}$ is called a Jordan algebra
if it satisfies 
\begin{center}
    $\big(x^2y\big)x\ =\ x^2\big(yx\big).$ 
\end{center}
A commutative algebra ${\rm A}$ is called a derived Jordan algebra
if the derived algebra ${\rm A}^{2}$ is a Jordan algebra.
\end{definition}

It is clear that any derived commutative associative algebra is a derived
Jordan algebra. The variety of derived Jordan algebras contains all Jordan
algebras.
Let ${\rm A}$ be a derived Jordan algebra. If either ${\rm A}^{2}=%
{\rm A}$ or ${\rm A}^{2}=0$, then ${\rm A}$\ is a Jordan algebra.

\begin{theoremA3}
    Let ${\mathcal J}$ be a complex $3$-dimensional derived Jordan algebra.
Then ${\mathcal J}$ is 
a Jordan algebra listed in {\rm Proposition \ref{3-dim Jord}}, or
a derived commutative associative algebra listed in {\rm Theorem A2},   
or it is isomorphic to one of the following algebras:

\begin{longtable}{lllllllllll}
 
${\mathcal J}_{01}^{\alpha}$ & $:$ & $e_{1}e_{1} = e_{1}$ & $e_{1}e_{2} = \frac{1}{2}e_{2}$ & $e_{2}e_{3} = e_{1}$ & $e_{3}e_{3} = e_{1} + \alpha e_{2}$ & \\ 
${\mathcal J}_{02}$ & $:$ & $e_{1}e_{1} = e_{1}$ & $e_{1}e_{2} = \frac{1}{2}e_{2}$ & $e_{2}e_{3} = e_{1}$ & $e_{3}e_{3} = e_{2}$ & \\
${\mathcal J}_{03}$ & $:$ & $e_{1}e_{1} = e_{1}$ & $e_{1}e_{2} = \frac{1}{2}e_{2}$ & $e_{2}e_{3} = e_{1}$ & & \\
${\mathcal J}_{04}^{\alpha}$ & $:$ & $e_{1}e_{1} = e_{1}$ & $e_{1}e_{2} = \frac{1}{2}e_{2}$ & $e_{2}e_{3} = e_{2}$ & $e_{3}e_{3} = \alpha e_{1}$ & \\
${\mathcal J}_{05}$ & $:$ & $e_{1}e_{1} = e_{1}$ & $e_{1}e_{2} = \frac{1}{2}e_{2}$ & $e_{2}e_{3} = e_{2}$ & $e_{3}e_{3} = -4e_{1} + e_{2}$ & \\
${\mathcal J}_{06}$ & $:$ & $e_{1}e_{1} = e_{1}$ & $e_{1}e_{2} = \frac{1}{2}e_{2}$ & $e_{3}e_{3} = e_{1}$ & & \\
${\mathcal J}_{07}$ & $:$ & $e_{1}e_{1} = e_{1}$ & $e_{1}e_{2} = \frac{1}{2}e_{2}$ & $e_{3}e_{3} = e_{2}$ & & \\
 
 \end{longtable}\noindent 
All algebras, excepting  ${\mathcal J}_{01}^{\alpha}\cong {\mathcal J}_{01}^{-
\alpha},$ are non-isomorphic.
\end{theoremA3}

\begin{proof}
 \underline{${\rm A}^{2}={\mathfrak J}_{04}$.} Then the nontrivial
multiplications in ${\rm A}$ is given as follows: 
\begin{longtable}{lcllll}
${\rm A}$ & $:$ & 
$e_{1}e_{1}=e_{1}$  & $e_{1}e_{3}=C_{131}e_{1}+C_{132}e_{2}$&$e_{2}e_{3}=C_{231}e_{1}+C_{232}e_{2}$\\
&& & $e_{1}e_{2}=\frac{1}{2}e_{2} $ &$e_{3}e_{3}=C_{331}e_{1}+C_{332}e_{2}$
\end{longtable}%
Then with the change of basis: 
$e_{1}^{\prime }=e_{1},$ \ 
$e_{2}^{\prime}=e_{2},$ \ 
$e_{3}^{\prime }=-C_{131}e_{1}-2C_{132}e_{2}+e_{3},$ 
we can make $C_{131}=C_{132}=0$. So we assume $C_{131}=C_{132}=0$. Hence the nontrivial
multiplications in ${\rm A}$ is given as follows: 
\begin{longtable}{lcllll}
${\rm A}$&$:$&$e_{1}e_{1}=e_{1}$&$e_{1}e_{2}=\frac{1}{2}e_{2}$&$e_{2}e_{3}=C_{231}e_{1}+C_{232}e_{2}$&$e_{3}e_{3}=C_{331}e_{1}+C_{332}e_{2}$
\end{longtable}%
Whence, ${\rm A}$ is a Jordan algebra if $%
C_{231}=C_{232}=C_{331}=C_{332}=0$.

\begin{itemize}
\item If $C_{231}\neq 0,$  then with the change of basis: 
\begin{center}
    $e_{1}^{\prime }=e_{1}+\frac{2}{3}C_{232} C^{-1}_{231}e_{2},$ \ 
    $e_{2}^{\prime }= {C^{-1}_{231}}e_{2},$ \ 
$e_{3}^{\prime }=-\frac{2}{3}C_{232}e_{1}-\frac{8}{9}C_{232}^{2}C^{-1}_{231}e_{2}+e_{3},$
\end{center} we can make $C_{231}=1$ and $C_{232}=0$. So we assume $C_{231}=1$ and $C_{232}=0,$  therefore the nontrivial multiplications in ${\rm A}$ is given as follows: 
\begin{longtable}{lcllllllll}
${\rm A}$&$:$&$e_{1}e_{1}=e_{1}$ &$e_{1}e_{2}=\frac{1}{2}e_{2}$ & $e_{2}e_{3}=e_{1}$& $e_{3}e_{3}=C_{331}e_{1}+C_{332}e_{2}.$
\end{longtable}

\begin{itemize}
\item If $C_{331}\neq 0,$ then, with the change of basis: 
$e_{1}^{\prime}=e_{1},$ \ 
$e_{2}^{\prime }={C^{\frac{1}{2}}_{331}}e_{2},$ \ 
${e_{3}^{\prime }={C^{-\frac{1}{2}}_{331}}e_{3}},$ we get the family of algebras ${\mathcal J}_{01}^{\alpha}$. Moreover, ${\mathcal J}_{01}^{\alpha}\cong {\mathcal J}_{01}^{\beta}$ if and only if $%
\alpha  =-\beta$.

\item If $C_{331}=0$ and $C_{332}\neq 0,$ then by changing 
$e_{1}^{\prime}=e_{1},$ \ 
$e_{2}^{\prime }=C^{\frac{1}{3}}_{332}e_{2},$ \ 
${e_{3}^{\prime }=C^{-\frac{1}{3}}_{332}e_{3}},$ we get the algebra ${\mathcal J}_{02}$.

\item If $C_{331}=C_{332}=0,$ we get the algebra ${\mathcal J}_{03}$.
\end{itemize}

\item If $C_{231}=0,$ then we have the following multiplication table:
\begin{longtable}{lcllllll}
${\rm A}$&$:$&$e_{1}e_{1}=e_{1}$&
$e_{1}e_{2}=\frac{1}{2}e_{2}$& 
$e_{2}e_{3}=C_{232}e_{2}$&
$e_{3}e_{3}=C_{331}e_{1}+C_{332}e_{2}$
\end{longtable}

\begin{itemize}
\item If $C_{232}\neq 0$ and $C_{331}+4C_{232}^{2}\neq 0,$ 
then with the change of basis: 
\begin{center}$
e_{1}^{\prime }=e_{1}+\frac{C_{332}}{C_{331}+4C_{232}^{2}}e_{2},$ \ 
$e_{2}^{\prime }=e_{2},$ \ 
$e_{3}^{\prime }=-\frac{2C_{332}}{C_{331}+4C_{232}^{2}}e_{2}+\frac{1}{C_{232}}e_{3}$\end{center} 
we get the family of algebras $%
{\mathcal J}_{04}^{\alpha \neq -4}$. Moreover, ${\mathcal J}_{04}^{\alpha}\cong {\mathcal J}_{04}^{\beta}$ if and only if $\alpha =\beta $.

\item If $C_{232}\neq 0,$
 $C_{331}+4C_{232}^{2}=0$ and $C_{332}=0$, then with the change of
basis: \begin{center} $e_{1}^{\prime }=e_{1},$ \ 
$e_{2}^{\prime }=e_{2},$ \ 
$e_{3}^{\prime }=C^{-1}_{232}e_{3},$
\end{center} we get the algebra ${\mathcal J}_{04}^{-4}$.

\item If $C_{232}\neq 0,$
 $C_{331}+4C_{232}^{2}=0$ and  $%
C_{332}\neq 0$, then with the change of basis: 
\begin{center}$e_{1}^{\prime}=e_{1},$ \ 
$e_{2}^{\prime }={C_{332}}{C_{232}^{-2}}e_{2},$ \ 
$e_{3}^{\prime }={C^{-1}_{232}}e_{3},$
\end{center} we get the algebra ${\mathcal J}_{05}$.

\item If $C_{232}=0$ and
 $C_{331}\neq 0,$ then with the change of basis: 
\begin{center}$e_{1}^{\prime }=e_{1}+{C_{332}}{C^{-1}_{331}}e_{2},$ \ 
$e_{2}^{\prime }=e_{2},$ \ 
${e_{3}^{\prime }={C^{-\frac{1}{2}}_{331}}e_{3}},$ 
\end{center}we get the algebra ${\mathcal J}_{06}$.

\item If $C_{232}=0$ and $C_{331}=0,$ 
then $C_{332}\neq 0$ and with the change of basis:
\begin{center}
$
e_{1}^{\prime }=e_{1},$\ 
$e_{2}^{\prime }=e_{2},$ \ 
$e_{3}^{\prime }={C^{-\frac{1}{2}}_{332}}e_{3},$ 
\end{center}we get the algebra ${\mathcal J}_{07}$.
\end{itemize}
\end{itemize}

\end{proof}

\subsection{Bicommutative algebras}\label{alg_bi}

\begin{definition}
An algebra is called a bicommutative algebra, if it satisfies the following identities%
\begin{longtable}{lcllcl}
$a(bc) $&$=$&$b(ac),$ & $  (ab)c $&$=$&$(ac)b.$
\end{longtable}
\end{definition}

\noindent It is easy to see if ${\rm A}$ is a bicommutative algebra, then%
\begin{eqnarray}
\big( \big( a\circ b\big) \circ c\big) \circ d &=&\big( \big( a\circ
b\big) \circ d\big) \circ c.  \label{4-assoc}
\end{eqnarray}%
Thus, 
${\rm A}^{+}$ is a derived commutative associative algebra. If ${\rm A}$ is commutative, then it is associative.
Let us call commutative algebras
satisfying the identity \eqref{4-assoc} as bicommutative$^+$ algebras.


\begin{proposition}
    
\label{B+}Let ${\rm A}$ be a complex $3$-dimensional bicommutative$^+$  algebra. Then ${\rm A}$
is a  commutative associative algebra listed in {\rm Proposition~\ref{3-dim assoc}} or
isomorphic to one of the following algebras:

\begin{longtable}{lllll}
${\rm M}_{04}^{\alpha}$ & $:$ & $e_{1}e_{3}=e_{1}$ & $e_{2}e_{3}=\alpha e_{2}$ \\ 
${\rm M}_{05}$ & $:$ & $e_{1}e_{3}=e_{1}$ & $e_{3}e_{3}=e_{2}$ \\ 
${\rm M}_{06}$ & $:$ & $e_{1}e_{3}=e_{1}+e_{2}$ & $e_{2}e_{3}=e_{2}$ \\ 
${\rm M}_{07}$ & $:$ & $e_{1}e_{3}=e_{1}$ & $e_{2}e_{2}=e_{1}$ \\ 
${\rm A}_{01}^{0,0,\gamma}$ & $:$ & $e_{1}e_{1}=e_{1}$ & $e_{2}e_{2}=e_{2}$ & $e_{3}e_{3}=e_{1}+\gamma e_{2}$ \\ 
${\rm A}_{03}^{\alpha}$ & $:$ & $e_{1}e_{1}=e_{1}$ & $e_{2}e_{3}=e_{2}$ & $e_{3}e_{3}=\alpha e_{1}$ \\ 
${\rm A}_{11}$ & $:$ & $e_{1}e_{1}=e_{1}$ & $e_{3}e_{3}=e_{1}+e_{2}$ \\ 
${\rm A}_{16}$ & $:$ & $e_{1}e_{1}=e_{1}$ & $e_{1}e_{2}=e_{2}$ & $e_{3}e_{3}=e_{1}+e_{2}$ \\ 
${\rm A}_{17}$ & $:$ & $e_{1}e_{1}=e_{1}$ & $e_{1}e_{2}=e_{2}$ & $e_{3}e_{3}=e_{1}$ \\ 
${\rm A}_{20}$ & $:$ & $e_{1}e_{1}=e_{1}$ & $e_{1}e_{2}=e_{2}$ & $e_{3}e_{3}=e_{2}$ \\ 
${\rm A}_{24}^{1}$ & $:$ & $e_{1}e_{1}=e_{2}$ & $e_{1}e_{3}=e_{1}$ & $e_{2}e_{3}=e_{2}$ \\ 
${\rm A}_{26}$ & $:$ & $e_{1}e_{1}=e_{2}$ & $e_{1}e_{3}=e_{1}+e_{2}$ & $e_{2}e_{3}=e_{2}$ \\ 
${\rm A}_{27}$ & $:$ & $e_{1}e_{1}=e_{2}$ & $e_{3}e_{3}=e_{1}$ \\ 
${\rm A}_{29}$ & $:$ & $e_{1}e_{1}=e_{1}$ & $e_{2}e_{3}=e_{1}$ \\ 
${\rm A}_{30}$ & $:$ & $e_{1}e_{1}=e_{1}$ & $e_{3}e_{3}=e_{1}$ \\ 
\end{longtable}
 
 \end{proposition}

\subsubsection{Preliminaries: the algebraic classification}\label{metbi}
In this subsection, we introduce the techniques used to obtain our main results 
$\big($the techniques are similar to those considered  in \cite{afm}$\big).$

\begin{definition}
Let $({\rm A},\cdot)$ be a bicommutative$^+$  algebra. 
Let ${\rm Z}^2({\rm A},{\rm A})$ be the set of all skew symmetric  bilinear maps $\theta\colon %
{\rm A}\times {\rm A} \to {\rm A}$ such that 
\begin{longtable}{ccccccccccccccccc}
        $(x\cdot y) \cdot z$&$+$& 
        $\theta(x,y)\cdot z$&$+$& 
        $\theta( x\cdot y,z)$&$+$& 
        $\theta(\theta(x,y),z)$&$=$ \\   
       && $(x \cdot z)\cdot y$&$+$& 
        $\theta(x,z)\cdot y$&$+$& 
        $\theta(x\cdot z,y)$&$+$& 
        $\theta(\theta(x,z),y),$\\
  
        $x \cdot (y\cdot z)$&$+$& 
        $x \cdot \theta(y,z)$&$+$& 
        $\theta(x, y\cdot z)$&$+$& 
        $\theta(x,\theta(y,z))$&$=$ \\   
       &&         $y \cdot (x\cdot z)$&$+$& 
        $y \cdot \theta(x,z)$&$+$& 
        $\theta(y, x\cdot z)$&$+$& 
        $\theta(y,\theta(x,z))$
        
\end{longtable}  
\end{definition}

\noindent  For $\theta \in {\rm Z}^2 (
{\rm A},{\rm A})$ we define on ${\rm A}$ a product $%
*_{\theta}\colon {\rm A}\times {\rm A}\to {\rm A}$ by 
$
x *_{\theta} y := \theta(x,y).  $

\begin{lemma}
Let $({\rm A},\cdot)$ be a bicommutative$^+$ algebra and $\theta \in {\rm Z}^2({\rm A},{\rm A})$. Then  $( {\rm A},\cdot_{\theta})$ is a
bicommutative algebra, 
where $x\cdot_{\theta} y := x \cdot y + x *_{\theta} y.$ 

\end{lemma}

Now, let $({\rm A},\cdot)$ be an algebra and $\mathrm{{Aut}(%
{\rm A})}$ be the automorphism group of ${\rm A}$ with respect to
product $\cdot$. Then $\mathrm{{Aut}({\rm A})}$ acts on ${\rm Z}^2({\rm A},%
{\rm A})$ by 
\begin{equation*}
(\theta \ast \phi)(x,y) := \phi^{-1}\bigl(\theta \bigl(\phi(x),\phi(y) \bigr)%
\bigr),
\end{equation*}
where $\phi \in \mathrm{{Aut}({\rm A})}$ and $\theta \in {\rm Z}^2({\rm A}
, {\rm A})$.

\begin{lemma}
 Let $({\rm A},\cdot)$ be a bicommutative$^+$ algebra and $\theta, \vartheta \in {\rm Z}^2 (%
{\rm A},{\rm A})$. Then the   algebras $({\rm A}, \cdot_{\theta})$ and $({\rm A}, \cdot_{\vartheta})$ are isomorphic if and only if there exists $\phi \in \mathrm{{Aut}({\rm A})}$ satisfying $\theta \ast \phi
=\vartheta $.
\end{lemma}

Hence, we have a procedure to classify the bicommutative algebras associated with a given bicommutative$^+$   algebra $({\rm A},\cdot)$. It consists of three steps:

\begin{enumerate}
\item[{\bf Step} $1.$] Compute ${\rm Z}^2 ({\rm A},{\rm A})$.

\item[{\bf Step} $2.$] Find the orbits of $\mathrm{{Aut}({\rm A})}$ on ${\rm Z}^2(
{\rm A},{\rm A})$.

\item[{\bf Step} $3.$] Take a representative $\theta$ for every orbit and  
construct an algebra $({\rm A}, \cdot_{\theta})$.
\end{enumerate}

Let us introduce the following notations. Let $\{e_1,\dots,e_n\}$ be a fixed basis of an algebra $({\rm A},\cdot)$. Define 
$\mathrm{\Lambda }^2({\rm A},\mathbb{C})$ to be the space of all 
  skew-symmetric bilinear forms on ${\rm A}$, that is, \begin{center}
 $\mathrm{\Lambda}^2(%
{\rm A},\mathbb{C}) := 
\big\langle \Delta_{ij} | 1\leq i < j \leq n \big\rangle$,
\end{center} where $\Delta_{ij}$ is the  skew-symmetric  bilinear form $%
\Delta_{ij}\colon {\rm A}\times {\rm A} \to \mathbb{C}$ defined by 
\begin{equation*}
\Delta_{ij}( e_l,e_m) :=\left\{ 
\begin{tabular}{rl}
$1,$ & if $(i,j) = (l,m),$   \\ 
$-1,$ & if    $(i,j) = (m,l),$   \\ 
$0,$ & otherwise.%
\end{tabular}
\right.
\end{equation*}
Now, if $\theta \in {\rm Z}^2 ({\rm A},{\rm A})$ then $\theta$ can be
uniquely written as $\theta (x,y) = \sum_{i=1}^n B_i(x,y)e_i$, where $%
B_1,\dots, B_n$ are   skew-symmetric  bilinear forms on ${\rm A}$. Also,
we may write $\theta = \big(B_{1},\dots ,B_n \big)$. Let $\phi^{-1}\in \mathrm{Aut}(%
{\rm A})$ be given by the matrix $( b_{ij})$. If $(\theta \ast
\phi)(x,y) = \sum_{i=1}^n B_i^{\prime }(x,y)e_i$, then 
\begin{center}
$B_i^{\prime }=
\sum_{j=1}^n b_{ij}\phi^t B_j\phi$, whenever $i \in \{1,\dots,n\}$.
\end{center}

\subsubsection{The classification Theorem}\label{TA2}
\begin{theoremA4}
Let ${\rm B}$ be a complex $3$-dimensional bicommutative algebra. Then $%
{\rm B}$ is a commutative associative algebra listed in Proposition \ref%
{3-dim assoc} or isomorphic to one of the following algebras\footnote{For receiving similar multiplication tables we have to apply the basis change $e_3:=\frac{1}{2}e_3$ in algebras ${\rm B}_{12}$--${\rm B}_{19}$, ${\rm B}_{22}$, ${\rm B}_{23}$.}:
\begin{longtable}{lllllll}
$\mathcal{G}_{00}$ & $:$ & $e_{2}e_{3}=e_{1}$ & $e_{3}e_{2}=-e_{1}$ \\
${\rm B}_{01}$ & $:$ & $e_{1}e_{1}=e_{1}$ & $e_{1}e_{2}=e_{2}$ & $e_{2}e_{1}=e_{2}$ \\
&& $e_{1}e_{3}=e_{2}$   & $e_{3}e_{1}=-e_{2}$ \\
${\rm B}_{02}$ & $:$ & $e_{1}e_{1}=e_{2}$ & $e_{1}e_{3}=e_{2}$ & $e_{3}e_{1}=-e_{2}$ \\
${\rm B}_{03}$ & $:$ & $e_{1}e_{1}=e_{2}$ & $e_{1}e_{2}=e_{3}$ & $e_{2}e_{1}=-e_{3}$ \\
${\rm B}_{04}^{\alpha \neq 0}$ & $:$ & $e_{1}e_{2}=\big( 1+\alpha \big) e_{3}$ & $e_{2}e_{1}=\big( 1-\alpha \big) e_{3}$ \\
${\rm B}_{05}^{\alpha \neq 0}$ & $:$ & $e_{1}e_{1}=e_{2}$ & $e_{1}e_{2}=\big( 1+\alpha \big) e_{3}$ & $e_{2}e_{1}=\big( 1-\alpha \big) e_{3}$ \\
${\rm B}_{06}^{\gamma}$ & $:$ & $e_{1}e_{1}=e_{1}$ & $e_{2}e_{2}=e_{2}$ & $e_{3}e_{3}=-e_{1}-\gamma^2 e_{2}$ & $e_{1}e_{3}=e_{1}$ \\
 &  & $e_{3}e_{1}=-e_{1}$ & $e_{2}e_{3}=\gamma e_{2}$ & $e_{3}e_{2}=-\gamma e_{2}$ \\


${\rm B}_{07}^{\gamma}$ & $:$ & $e_{1}e_{1}=e_{1}$ & $e_{1}e_{3}=\gamma e_{1}$ & $e_{3}e_{1}=-\gamma e_{1}$ \\
&& $e_{2}e_{3}=2e_{2}$   & $e_{3}e_{3}=-\gamma^2 e_{1}$ \\
${\rm B}_{08}^{\gamma}$ & $:$ & $e_{1}e_{1}=e_{1}$ & $e_{1}e_{3}=\gamma e_{1}$ & $e_{3}e_{1}=-\gamma e_{1}$ \\
&& $e_{3}e_{2}=2e_{2}$   & $e_{3}e_{3}=-\gamma^2 e_{1}$ \\


${\rm B}_{09}$ & $:$ & $e_{1}e_{1}=e_{1}$ & $e_{3}e_{3}=-e_{1}-e_{2}$ & $e_{1}e_{3}=e_{1}$ & $e_{3}e_{1}=-e_{1}$ \\

${\rm B}_{10}$ & $:$ & $e_{1}e_{1}=e_{1}$ & $e_{1}e_{2}=e_{2}$ & $e_{2}e_{1}=e_{2}$ & $e_{3}e_{3}=-e_{1}-e_{2}$ \\
 &  & $e_{1}e_{3}=e_{1}+\frac{1}{2}e_{2}$
& $e_{3}e_{1}=-e_{1}-\frac{1}{2} e_{2}$ & $e_{2}e_{3}=e_{2}$ & $e_{3}e_{2}=-e_{2}$ \\

${\rm B}_{11}$ & $:$ & $e_{1}e_{1}=e_{1}$ & $e_{1}e_{2}=e_{2}$& $e_{2}e_{1}=e_{2}$ & $e_{3}e_{3}=-e_{1}$ \\
 &  &  $e_{1}e_{3}=e_{1}$  & $e_{3}e_{1}=-e_{1}$ & $e_{2}e_{3}=e_{2}$ & $e_{3}e_{2}=-e_{2}$ \\

${\rm B}_{12}$ & $:$ & $e_{1}e_{1}=e_{2}$ & $e_{1}e_{3}=e_{1}$ & $e_{2}e_{3}=e_{2}$ \\
${\rm B}_{13}$ & $:$ & $e_{1}e_{1}=e_{2}$ & $e_{3}e_{1}=e_{1}$ & $e_{3}e_{2}=e_{2}$ \\
${\rm B}_{14}$ & $:$ & $e_{1}e_{1}=e_{2}$ & $e_{1}e_{3}=e_{1}+e_{2}$ & $e_{2}e_{3}=e_{2}$ \\
${\rm B}_{15}$ & $:$ & $e_{1}e_{1}=e_{2}$ & $e_{3}e_{1}=e_{1}+e_{2}$ & $e_{3}e_{2}=e_{2}$ \\
${\rm B}_{16}^{\alpha}$ & $:$ & $e_{1}e_{3}=e_{1}$ & $e_{2}e_{3}=\alpha e_{2}$ \\
${\rm B}_{17}^{\alpha \neq 0}$ & $:$ & $e_{1}e_{3}=e_{1}$ & $e_{3}e_{2}=\alpha e_{2}$ \\
${\rm B}_{18}^{\alpha \neq \pm 1}$ & $:$ & $e_{3}e_{1}=e_{1}$ & $e_{2}e_{3}=\alpha e_{2}$ \\
${\rm B}_{19}^{\alpha \neq 0}$ & $:$ & $e_{3}e_{1}=e_{1}$ & $e_{3}e_{2}=\alpha e_{2}$ \\
${\rm B}_{20}$ & $:$ & $e_{1}e_{3}=2e_{1}$ & $e_{3}e_{3}=e_{2}$ \\
${\rm B}_{21}$ & $:$ & $e_{3}e_{1}=2e_{1}$ & $e_{3}e_{3}=e_{2}$ \\
${\rm B}_{22}$ & $:$ & $e_{1}e_{3}=e_{1}+e_{2}$ & $e_{2}e_{3}=e_{2}$ \\
${\rm B}_{23}$ & $:$ & $e_{3}e_{1}=e_{1}+e_{2}$ & $e_{3}e_{2}=e_{2}$ \\
${\rm B}_{24}$ & $:$ & $e_{1}e_{1}=e_{1}$ & $e_{1}e_{3}=e_{1}$ & $e_{3}e_{1}=-e_{1}$ & $e_{3}e_{3}=-e_{1}$ \\
\end{longtable}
\noindent
All listed algebras are non-isomorphic, excepting:
\begin{center}
$
{\rm B}_{04}^{\alpha}\cong {\rm B}_{04}^{-\alpha}, \
{\rm B}_{06}^{\alpha} \cong {\rm B}_{06}^{\alpha^{-1}}, \ 
{\rm B}_{16}^{\alpha} \cong {\rm B}_{16}^{\alpha^{-1}}, \  
{\rm B}_{19}^{\alpha} \cong {\rm B}_{19}^{\alpha^{-1}}$.
\end{center}
\end{theoremA4}
\subsubsection{Proof of  Theorem A4}

If ${\rm B}^{+}$ is trivial, then ${\rm B}$\ is a metabelian Lie algebra, and we obtain the algebra $\mathcal{G}_{00}$. If ${\rm B}^{+}$ is nontrivial, then ${\rm B}^{+}$ is isomorphic to one of the algebras
listed in Proposition~\ref{B+}. Now, assume that ${\rm B}^{+}$ is one of the algebras
in Proposition~\ref{B+}. By direct computation of the set
${\rm Z}^{2}({\rm B}^{+}, {\rm B}^{+})$, we have the following observations:

\begin{itemize}
    \item If ${\rm B}^{+}$ is associative, then
    \[
        {\rm Z}^{2}({\rm B}^{+}, {\rm B}^{+}) = \{0\} \quad \text{for} \quad
        {\rm B}^{+} \in \big \{ {\rm J}_{01},\ {\rm J}_{03},\  {\rm J}_{07},\ {\rm J}_{08},\ {\rm J}_{09},\ {\rm J}_{10},\ {\rm J}_{11} \big\}.
    \]
   
    \item If ${\rm B}^{+}$ is non-associative, then
    \[
        {\rm Z}^{2}({\rm B}^{+}, {\rm B}^{+}) = \varnothing \quad \text{for} \quad
        {\rm B}^{+} \in  \big\{ {\rm A}_{20},\ {\rm A}_{27},\ {\rm M}_{07},\ {\rm A}_{29} \big\}.
    \]
   
\end{itemize}

Therefore, we focus on the following cases:

\begin{enumerate}[I.]
    \item 
\underline{${\rm B}^{+}={\rm J}_{02}$.}\ Let 
$0 \neq \theta =\big(
B_{1},B_{2},B_{3}\big) \in {\rm Z}^{2}\big( 
{\rm B}^{+},{\rm B}^{+}\big) $. Then 
$\theta =\big( 0,\ \alpha \Delta _{13},\ 0\big) $ for some $\alpha \in \mathbb{C}^{\ast}.$ The
automorphism group of ${\rm J}_{02}$ consists of the automorphisms $\phi $ given by a matrix of
the following form:%
\begin{equation*}
\phi =%
\begin{pmatrix}
1 & 0 & 0 \\ 
0 & a_{22} & 0 \\ 
0 & 0 & a_{33}%
\end{pmatrix}%
.
\end{equation*}%
Let $\phi =\bigl(a_{ij}\bigr)\in $ $\text{Aut}\big( {\rm J}_{02}\big) 
$. Then $\theta \ast \phi =\big( 0,\ \beta \Delta _{13},\ 0\big),$ 
where $\beta =\alpha a_{33}a^{-1}_{22}$. 
Let $\phi $\ be the following
automorphism: 
\begin{equation*}
\phi =%
\begin{pmatrix}
1 & 0 & 0 \\ 
0 & \alpha & 0 \\ 
0 & 0 & 1%
\end{pmatrix}%
.
\end{equation*}%
Then $\theta \ast \phi =\big( 0,\ \Delta _{13},\ 0\big) $. So we get the
algebra ${\rm B}_{01}$.

\bigskip

\item 
\underline{$\rm{B}^{+}={\rm J}_{04}$.}\ Let 
$0\neq \theta =\big(
B_{1},B_{2},B_{3}\big) \in {\rm Z}^{2}\big( 
\rm{B}^{+},\rm{B}^{+}\big) $. Then $\theta \in \left\{ \eta
_{1},\eta _{2},\eta _{3}\right\},$ where
\begin{longtable}{lcl}
$\eta _{1} $&$=$&$\big( 0,\alpha _{1}\Delta _{13},0\big) ,$ \\
$\eta _{2} $&$=$&$\big( 0,0,\alpha _{1}\Delta _{12}\big) ,$ \\
$\eta _{3} $&$=$&$\big( 0,\alpha _{1}\Delta _{12}+\alpha _{2}\Delta _{13},-\frac{%
\alpha _{1}^{2}}{\alpha _{2}}\Delta _{12}-\alpha _{1}\Delta _{13}\big)
_{\alpha _{1}\alpha _{2}\neq 0},$
\end{longtable}%
for some $\alpha _{1},\alpha _{2}\in \mathbb{C}$. The automorphism group of $%
{\rm J}_{04}$  consists of the automorphisms $\phi $ given by a matrix of the following form:%
\begin{equation*}
\phi =%
\begin{pmatrix}
a_{11} & 0 & 0 \\ 
a_{21} & a_{11}^{2} & a_{23} \\ 
a_{31} & 0 & a_{33}%
\end{pmatrix}%
.
\end{equation*}

\begin{itemize}
\item $\theta =\eta _{1}$. Let $\phi $\ be the following automorphism: 
\begin{equation*}
\phi =%
\begin{pmatrix}
\alpha _{1} & 0 & 0 \\ 
0 & \alpha _{1}^{2} & 0 \\ 
0 & 0 & 1%
\end{pmatrix}%
.
\end{equation*}%
Then $\theta \ast \phi =\big( 0,\ \Delta _{13},\ 0\big) $. So we get the
algebra ${\rm B}_{02}$.

\item $\theta =\eta _{2}$. Let $\phi $\ be the following automorphism: 
\begin{equation*}
\phi =%
\begin{pmatrix}
1 & 0 & 0 \\ 
0 & 1 & 0 \\ 
0 & 0 & \alpha _{1}%
\end{pmatrix}%
.
\end{equation*}%
Then $\theta \ast \phi =\big( 0,\ 0,\ \Delta _{12}\big) $. So we get the
algebra ${\rm B}_{03}$.

\item $\theta =\eta _{3}$. Let $\phi $\ be the following automorphism:%
\begin{equation*}
\begin{pmatrix}
1 & 0 & 0 \\ 
0 & 1 & \alpha _{1} \\ 
0 & 0 & -{\alpha _{1}^{2}}{\alpha^{-1}_{2}}%
\end{pmatrix}%
\end{equation*}%
Then $\theta \ast \phi =\big( 0,\ 0,\ \Delta _{12}\big) $. So we get the
algebra ${\rm B}_{03}$.
\end{itemize}

\item 
\underline{${\rm B}^{+}={\rm J}_{05}$.}\ 
Let $0 \neq \theta =\big(
B_{1},B_{2},B_{3}\big) \in {\rm Z}^{2}\big( 
{\rm B}^{+},{\rm B}^{+}\big) $. 
Then $\theta =\big( 0,\ 0,\ \alpha\Delta _{12}\big) $ for some $\alpha \in \mathbb{C}^{\ast }$. The
automorphism group of ${\rm J}_{05}$ consists of the automorphisms $\phi $ given by a matrix of
the following form:%
\begin{equation*}
\phi =%
\begin{pmatrix}
\varepsilon a_{11} & \epsilon a_{12} & 0 \\ 
\epsilon a_{21} & \varepsilon a_{22} & 0 \\ 
a_{31} & a_{32} & \varepsilon a_{11}a_{22}+\epsilon a_{12}a_{21}%
\end{pmatrix}%
: \ \big( \varepsilon ,\epsilon \big) \in \big\{\big( 1,0\big) ,\big(0,1\big) \big\}.
\end{equation*}%
Let $\phi =\bigl(a_{ij}\bigr)\in $ $\text{Aut}\big( {\rm J}_{05}\big) 
$. Then $\theta \ast \phi =\big( 0,\ 0,\ \beta \Delta _{12}\big) $ where $%
\beta ^{2}=\alpha ^{2}$. So we get the family of algebras ${\rm B}_{04}^{\alpha
\neq 0}$. Moreover, ${\rm B}_{04}^{\alpha}\cong {\rm B}_{04}^{\beta}$ if and only if $\alpha =\pm\beta$.

\item 
\underline{${\rm B}^{+}={\rm J}_{06}$.}\ Let
$0 \neq \theta =\big(
B_{1},B_{2},B_{3}\big) \in {\rm Z}^{2}\big( 
{\rm B}^{+},{\rm B}^{+}\big) $. Then $\theta =\big( 0,\ 0,\ \alpha
\Delta _{12}\big) $ for some $\alpha \in \mathbb{C}^{\ast }$. The
automorphism group of ${\rm J}_{06}$, $\text{Aut}\big( {\rm J}%
_{06}\big) $, consists of the automorphisms $\phi $ given by a matrix of
the following form:%
\begin{equation*}
\phi =%
\begin{pmatrix}
a_{11} & 0 & 0 \\ 
a_{21} & a_{11}^{2} & 0 \\ 
a_{31} & 2a_{11}a_{21} & a_{11}^{3}%
\end{pmatrix}%
.
\end{equation*}%
Let $\phi =\bigl(a_{ij}\bigr)\in $ $\text{Aut}\big( {\rm J}_{06}\big) 
$. Then $\theta \ast \phi =\theta $. So we get the family of algebras ${\rm B}%
_{05}^{\alpha \neq 0}$. Moreover, ${\rm B}_{05}^{\alpha}\cong {\rm B}_{05}^{\beta}$ if and only if $\alpha =\beta$.

\item \underline{$\mathrm{B}^{+}=\mathrm{A}_{01}^{0,0,\gamma}$.}
Consider the change of basis $e_{1}^{\prime }=e_{1},$\ 
$e_{2}^{\prime }=e_{2},$\ $e_{3}^{\prime }={\rm i}e_{3}$. Then%
\begin{equation*}
\mathrm{A}_{01}^{0,0,\gamma }\cong \mathrm{A}_{\gamma }^{\prime }\ :\
e_{1}^{\prime }e_{1}^{\prime }=e_{1}^{\prime },\ e_{2}^{\prime
}e_{2}^{\prime }=e_{2}^{\prime },\ e_{3}^{\prime }e_{3}^{\prime
}=-e_{1}^{\prime }-\gamma e_{2}^{\prime }.
\end{equation*}%

Hence, we may assume $\mathrm{B}^{+}=\mathrm{A}_{\gamma}'$.
Let $\theta =%
\big(B_{1},B_{2},B_{3}\big)$ be an arbitrary element of $\mathrm{Z}^{2}\big(%
\mathrm{B}^{+},\mathrm{B}^{+}\big)$. Then $\theta =\big(\epsilon \Delta
_{13},\alpha \Delta _{23},0\big)$
where $\epsilon^{2}=1$ and $\alpha^{2}=\gamma$.

Therefore, the bicommutative algebra structures defined on
$\mathrm{A}_{\gamma}'$ are given by the following families:
\begin{longtable}{lclcllcllcllcl}
$X^{\alpha } $&$:$& $
e_{1}^{\prime }e_{1}^{\prime }$&$=$&$e_{1}^{\prime },$&$
e_{2}^{\prime }e_{2}^{\prime }$&$=$&$e_{2}^{\prime },$&$ 
e_{3}^{\prime}e_{3}^{\prime }$&$=$&$-e_{1}^{\prime }-\alpha ^{2}e_{2}^{\prime },$&$ 
e_{1}^{\prime}e_{3}^{\prime }$&$=$&$e_{1}^{\prime },$\\&& $
e_{3}^{\prime }e_{1}^{\prime}$&$=$&$-e_{1}^{\prime },$&$
e_{2}^{\prime }e_{3}^{\prime }$&$=$&$\alpha e_{2}^{\prime },$&$
e_{3}^{\prime }e_{2}^{\prime }$&$=$&$-\alpha e_{2}^{\prime }.$ \\

$Y^{\alpha }$ &$:$& $
e_{1}^{\prime }e_{1}^{\prime }$&$=$&$e_{1}^{\prime },$&$
e_{2}^{\prime }e_{2}^{\prime }$&$=$&$e_{2}^{\prime },$&$ 
e_{3}^{\prime}e_{3}^{\prime }$&$=$&$-e_{1}^{\prime }-\alpha ^{2}e_{2}^{\prime },$&$ 
e_{1}^{\prime}e_{3}^{\prime }$&$=$&$-e_{1}^{\prime },$\\&& $ 
e_{3}^{\prime }e_{1}^{\prime}$&$=$&$e_{1}^{\prime },$&$
e_{2}^{\prime }e_{3}^{\prime }$&$=$&$\alpha e_{2}^{\prime },$&$
e_{3}^{\prime }e_{2}^{\prime }$&$=$&$-\alpha e_{2}^{\prime }.$
\end{longtable}
\noindent Moreover,
$X^{\alpha} \cong X^{\alpha^{-1}}
\cong Y^{-\alpha}.$
Hence, this case yields the family of algebras $\mathrm{B}_{06}^{\gamma}$.

\item \underline{$\mathrm{B}^{+} = \mathrm{A}_{03}^{-\gamma}$.} 
Let $0 \neq \theta = (B_1, B_2, B_3) \in \mathrm{Z}^2(\mathrm{B}^{+}, \mathrm{B}^{+})$. 
Then 
$\theta =\big(\alpha \Delta _{13},\ \epsilon\Delta _{23},\ 0\big)$ with $\alpha^2 = \gamma$  and  $\epsilon^2 = 1$. Hence, the bicommutative structures defined on $\mathrm{A}_{03}^{-\gamma}$ are given by the
following families:
\begin{align*}
X^{\alpha} &: \quad 
e_1 e_1 = e_1, \quad e_1 e_3 = \alpha e_1, \quad e_3 e_1 = -\alpha e_1, \quad 
e_2 e_3 = 2 e_2, \quad e_3 e_3 = -\alpha^2 e_1, \\
Y^{\alpha} &: \quad 
e_1 e_1 = e_1, \quad e_1 e_3 = \alpha e_1, \quad e_3 e_1 = -\alpha e_1, \quad 
e_3 e_2 = 2 e_2, \quad e_3 e_3 = -\alpha^2 e_1.
\end{align*}
Then $X^{\alpha}$ and $Y^{\beta}$ are not isomorphic. Moreover, 
$X^{\alpha} \cong X^{\beta}$ if and only if $\alpha = \beta$, and 
$Y^{\alpha} \cong Y^{\beta}$ if and only if $\alpha = \beta$. We obtain the families of algebras $\mathrm{B}_{07}^{\gamma}$ 
and $\mathrm{B}_{08}^{\gamma}$.

\item 
\underline{${\rm B}^{+}={\rm A}_{11}$.} Consider the change of basis 
$e_{1}^{\prime }=e_{1},$ 
$e_{2}^{\prime }=e_{2},$ 
$e_{3}^{\prime }={\rm i}e_{3}$. Then%
\begin{equation*}
{\rm A}_{11}\cong {\rm A}_{11}^{\prime }\ :\ 
e_{1}^{\prime}e_{1}^{\prime }=e_{1}^{\prime },\ 
e_{3}^{\prime }e_{3}^{\prime}=-e_{1}^{\prime }-e_{2}^{\prime }.
\end{equation*}%
So we assume ${\rm B}^{+}={\rm A}_{11}^{\prime }$. Let 
$0\neq \theta
=\big( B_{1},B_{2},B_{3}\big) \in {\rm Z}^{2}\big( 
{\rm B}^{+},{\rm B}^{+}\big) $. Then 
$\theta =\big( \alpha \Delta_{13},\ 0,\ 0\big)$ 
such that $\alpha ^{2}=1$. The automorphism group of $%
{\rm A}_{11}^{\prime }$  consists of the automorphisms $\phi $ given by a matrix of the
following form:%
\begin{equation*}
\phi =%
\begin{pmatrix}
1 & 0 & 0 \\ 
0 & 1 & a_{23} \\ 
0 & 0 & \epsilon%
\end{pmatrix}%
:\ \epsilon ^{2}=1.
\end{equation*}%
Let $\phi =\bigl(a_{ij}\bigr)\in $ $\text{Aut}\big( {\rm A}
_{11}^{\prime }\big) $. Then $\theta \ast \phi =\big( \epsilon \alpha
\Delta _{13},\ 0,\ 0\big) $. 
So we may choose $\epsilon $ such that $\epsilon
\alpha =1$ and therefore we obtain the algebra ${\rm B}_{09}$.

\item 
\underline{${\rm B}^{+}={\rm A}_{16}$.} Consider the change of basis 
$e_{1}^{\prime }=e_{1},$ 
$e_{2}^{\prime }=e_{2},$ 
$e_{3}^{\prime }={\rm i}e_{3}$. Then%
\begin{equation*}
{\rm A}_{16}\cong {\rm A}_{16}^{\prime }\ :\ 
e_{1}^{\prime}e_{1}^{\prime }=e_{1}^{\prime },\ 
e_{1}^{\prime }e_{2}^{\prime}=e_{2}^{\prime },\ 
e_{3}^{\prime }e_{3}^{\prime }=-e_{1}^{\prime
}-e_{2}^{\prime }.
\end{equation*}%
So we assume ${\rm B}^{+}={\rm A}_{16}^{\prime }$. Let 
$0\neq \theta =\big( B_{1},B_{2},B_{3}\big) \in {\rm Z}^{2}\big( 
{\rm B}^{+},{\rm B}^{+}\big) $. Then $\theta =\big( \alpha \Delta
_{13},\ \frac{1}{2}\alpha \Delta _{13}+\alpha \Delta _{23},\ 0\big) $ such
that $\alpha ^{2}=1$. The automorphism group of ${\rm A}_{16}^{\prime }$ 
 consists of the
automorphisms $\phi $ given by a matrix of the following form:%
\begin{equation*}
\phi _{1}=%
\begin{pmatrix}
1 & 0 & 0 \\ 
0 & 1 & 0 \\ 
0 & 0 & 1%
\end{pmatrix}%
\allowbreak,\ 
\phi _{2}=%
\begin{pmatrix}
1 & 0 & 0 \\ 
0 & 1 & 0 \\ 
0 & 0 & -1%
\end{pmatrix}%
.
\end{equation*}%
Then $\big( \Delta _{13},\ \frac{1}{2}\Delta _{13}+\Delta _{23},\ 0\big) \ast
\phi _{2}=\big( -\Delta _{13},\ -\frac{1}{2}\Delta _{13}-\Delta
_{23},\ 0\big) $ and we have   ${\rm B}_{10}$.

\item 
\underline{${\rm B}^{+}={\rm A}_{17}$.} Consider the change of basis 
$e_{1}^{\prime }=e_{1},$ 
$e_{2}^{\prime }=e_{2},$ 
$e_{3}^{\prime }={\rm i}e_{3}$. Then%
\begin{equation*}
{\rm A}_{17}\cong {\rm A}_{17}^{\prime }\ :\ 
e_{1}^{\prime}e_{1}^{\prime }=e_{1}^{\prime },\ 
e_{1}^{\prime }e_{2}^{\prime}=e_{2}^{\prime },\ 
e_{3}^{\prime }e_{3}^{\prime }=-e_{1}^{\prime }.
\end{equation*}%
So we assume ${\rm B}^{+}={\rm A}_{17}^{\prime }$. Let 
$0\neq \theta =\big( B_{1},B_{2},B_{3}\big) \in {\rm Z}^{2}\big( 
{\rm B}^{+},{\rm B}^{+}\big) $. Then $\theta =\big( \alpha \Delta_{13},\ \alpha \Delta _{23},\ 0\big) $ such that $\alpha ^{2}=1$. The
automorphism group of ${\rm A}_{17}^{\prime }$ consists of the automorphisms $\phi $
given by a matrix of the following forms:%
\begin{equation*}
\phi _{1} =%
\begin{pmatrix}
1 & 0 & 0 \\ 
0 & a_{22} & 0 \\ 
0 & 0 &  1%
\end{pmatrix}%
, \ 
\phi _{2} =%
\begin{pmatrix}
1 & 0 & 0 \\ 
0 & a_{22} & 0 \\ 
0 & 0 & - 1%
\end{pmatrix}%
.
\end{equation*}%
Then $\big(  \Delta
_{13}, \Delta _{23},0\big) \ast \phi _{2} =-\big(  \Delta
_{13}, \Delta _{23},0\big) $. So we get the algebra ${\rm B}_{11}$.

\item 
\underline{${\rm B}^{+}={\rm A}_{24}^{1}$.} Let 
$0 \neq \theta =\big(
B_{1},B_{2},B_{3}\big) \in {\rm Z}^{2}\big( {\rm B
}^{+},{\rm B}^{+}\big) $. Then $\theta =\big( \alpha \Delta
_{13},\alpha \Delta _{23},0\big) $ such that $\alpha ^{2}=1$. The
automorphism group of ${\rm A}_{24}^{1}$  consists of the automorphisms $\phi $ given by a matrix
of the following form:%
\begin{equation*}
\phi =%
\begin{pmatrix}
a_{11} & 0 & 0 \\ 
a_{21} & a_{11}^{2} & 0 \\ 
0 & 0 & 1%
\end{pmatrix}%
.
\end{equation*}%
Then $\theta \ast \phi =\theta $ for any $\phi \in \text{Aut}\big( {\rm  
A}_{24}^{1}\big) $. So we get the algebras ${\rm B}_{12}$ and ${\rm B}_{13}$.

\item 
\underline{${\rm B}^{+}={\rm A}_{26}$.} Let $0\neq \theta =\big(
B_{1},B_{2},B_{3}\big) \in {\rm Z}^{2}\big( {\rm B
}^{+},{\rm B}^{+}\big) $. Then 
$\theta =\big( \alpha \Delta_{13}, \ \alpha \Delta _{13}+\alpha \Delta _{23},\ 0\big) $ such that $\alpha
^{2}=1$. The automorphism group of ${\rm A}_{26}$ consists of the automorphisms $\phi $ given by a
matrix of the following form:%
\begin{equation*}
\phi =%
\begin{pmatrix}
1 & 0 & 0 \\ 
a_{21} & 1 & 0 \\ 
0 & 0 & 1%
\end{pmatrix}%
.
\end{equation*}%
Then $\theta \ast \phi =\theta $ for any $\phi \in \text{Aut}\big( {\rm  
A}_{26}\big) $. So we get the algebras ${\rm B}_{14}$ and ${\rm B}_{15}$.

\item 
\underline{${\rm B}^{+}={\rm M}_{04}^{\alpha \notin \{ 0,\pm 1\}}$.} 
Let 
$0\neq \theta =\big( B_{1},B_{2},B_{3}\big) \in 
{\rm Z}^{2}\big( {\rm B}^{+},{\rm B}^{+}\big) $. Then $\theta =\big(
\alpha _{1}\Delta _{13},\ \alpha _{2}\Delta _{23},\ 0\big) $ such that $\alpha
_{1}^{2}=1$ and $\alpha _{2}^{2}=\alpha ^{2}$. The automorphism group of ${\rm  
M}_{04}^{\alpha \notin \{ 0,\pm 1\}}$  consists of the automorphisms $\phi $ given by a
matrix of the following form:%
\begin{equation*}
\phi =%
\begin{pmatrix}
a_{11} & 0 & 0 \\ 
0 & a_{22} & 0 \\ 
0 & 0 & 1%
\end{pmatrix}%
.
\end{equation*}%
Then $\theta \ast \phi =\theta $ for any $\phi \in \text{Aut}\big( {\rm  
M}_{04}^{\alpha \notin\{ 0,\pm 1\}}\big) $. 
So we get the families of algebras ${\rm B}_{16}^{\alpha \notin\{ 0,\pm 1\}},$ 
${\rm B}_{17}^{\alpha \notin \{ 0,\pm 1\}},$ 
${\rm B}_{18}^{\alpha \notin \{ 0,\pm 1\}},$ and 
${\rm B}_{19}^{\alpha \notin \{ 0,\pm 1\}}$. Moreover, ${\rm B}_{16}^{\alpha}\cong {\rm B}_{16}^{\beta}$ if and only if $\alpha =\beta$ or $\alpha =\beta^{-1}$; ${\rm B}_{17}^{\alpha}\cong {\rm B}_{17}^{\beta}$ if and only if $\alpha =\beta$; ${\rm B}_{18}^{\alpha}\cong {\rm B}_{18}^{\beta}$ if and only if $\alpha =\beta$; ${\rm B}_{19}^{\alpha}\cong {\rm B}_{19}^{\beta}$ if and only if $\alpha =\beta$ or $\alpha =\beta^{-1}$.

\item 
\underline{${\rm B}^{+}={\rm M}_{04}^{-1}$.} Let 
$0 \neq \theta
=\big( B_{1},B_{2},B_{3}\big) \in {\rm Z}^{2}\big( 
{\rm B}^{+},{\rm B}^{+}\big) $. Then 
$\theta =\big(\alpha_{1}\Delta _{13},\ \alpha _{2}\Delta _{23},\ 0\big) $ such that $\alpha
_{1}^{2}=1$ and $\alpha _{2}^{2}=1$. The automorphism group of ${\rm M}_{04}^{-1}$  consists of the automorphisms $\phi $ given by a matrix of the following
forms:%
\begin{equation*}
\phi _{1}=%
\begin{pmatrix}
a_{11} & 0 & 0 \\ 
0 & a_{22} & 0 \\ 
0 & 0 & 1%
\end{pmatrix}%
,\ 
\phi _{2}=%
\begin{pmatrix}
0 & a_{12} & 0 \\ 
a_{21} & 0 & 0 \\ 
0 & 0 & -1%
\end{pmatrix}%
.
\end{equation*}%
Then $\big( \Delta _{13},\ \Delta _{23},\ 0\big)\ast \phi _{2} =
\big(-\Delta _{13},\ -\Delta _{23},\ 0\big) $ and 
 we have  ${\rm B}_{16}^{-1},$ \ 
${\rm B}_{17}^{-1},$ and 
${\rm B}_{19}^{-1}$.

\item 
\underline{${\rm B}^{+}={\rm M}_{04}^{0}$.} Let 
$0\neq \theta
=\big( B_{1},B_{2},B_{3}\big) \in {\rm Z}^{2}\big( 
{\rm B}^{+},{\rm B}^{+}\big) $. Then $\theta =\big( \alpha \Delta
_{13},\ 0,\ 0\big) $ such that $\alpha ^{2}=1$. The automorphism group of $%
{\rm M}_{04}^{0}$  consists of the automorphisms $\phi $ given by a matrix of the
following form:%
\begin{equation*}
\phi =%
\begin{pmatrix}
a_{11} & 0 & 0 \\ 
0 & a_{22} & a_{23} \\ 
0 & 0 & 1%
\end{pmatrix}%
.
\end{equation*}%
Then $\theta \ast \phi =\theta $ for any $\phi \in \text{Aut}\big( {\rm  
M}_{04}^{0}\big) $. So we get the algebras ${\rm B}_{16}^{0}$ and ${\rm B}_{18}^{0}$.

\item
\underline{${\rm B}^{+}={\rm M}_{04}^{1}$.} Let 
$0\neq \theta=\big( B_{1},B_{2},B_{3}\big) \in {\rm Z}^{2}\big( 
{\rm B}^{+},{\rm B}^{+}\big) $. Then $\theta \in \left\{ \eta
_{1},\ldots ,\eta _{5}\right\} $ where%
\begin{longtable}{lcl}
$\eta _{1} $&$=$&$\big( \alpha _{1}\Delta _{13}+(1-\alpha _{1}^{2})\alpha^{-1}
_{2}\Delta _{23},\ \alpha _{2}\Delta _{13}-\alpha _{1}\Delta _{23},\ 0\big)
_{\alpha _{2}\neq 0},$ \\
$\eta _{2} $&$=$&$\big( \Delta _{13}+\alpha _{1}\Delta _{23},\ 
-\Delta_{23}, \ 0\big) ,$ \\
$\eta _{3} $&$=$&$\big( -\Delta _{13}+\alpha _{1}\Delta _{23},\ 
\Delta_{23},\ 0\big),$ \\
$\eta _{4} $&$=$&$\big( \Delta _{13},\ \Delta _{23},\ 0\big) ,$ \\
$\eta _{5} $&$=$&$\big( -\Delta _{13},\ -\Delta _{23},\ 0\big) .$
\end{longtable}%
The automorphism group of ${\rm M}_{04}^{1}$ consists of the automorphisms $\phi $
given by a matrix of the following form:%
\begin{equation*}
\phi =%
\begin{pmatrix}
a_{11} & a_{12} & 0 \\ 
a_{21} & a_{22} & 0 \\ 
0 & 0 & 1%
\end{pmatrix}%
.
\end{equation*}%
Let $\eta =\big( \gamma _{1}\Delta _{13}+\gamma _{2}\Delta _{23},\ 
\gamma_{3}\Delta _{13}+\gamma _{4}\Delta _{23},\ 0\big) $ and 
$\phi =\bigl(a_{ij}%
\bigr)\in $ $\text{Aut}\big( {\rm M}_{04}^{1}\big) $. Then 
\begin{center}$%
\eta \ast \phi =\big( \gamma _{1}^{\prime }\Delta _{13}+\gamma _{2}^{\prime
}\Delta _{23},\ \gamma _{3}^{\prime }\Delta _{13}+\gamma _{4}^{\prime }\Delta
_{23},\ 0\big) $\end{center}
where 
\(
M_{\eta\ast\phi}=\phi^{-1}M_{\eta}\phi
\),
with 
\(
M_{\eta}=\begin{pmatrix}\gamma_{1}&\gamma_{2}\\[2pt]\gamma_{3}&\gamma_{4}\end{pmatrix}
\)
and
\(
M_{\eta\ast\phi}=\begin{pmatrix}\gamma_{1}'&\gamma_{2}'\\[2pt]\gamma_{3}'&\gamma_{4}'\end{pmatrix}.
\)

Hence we may assume $M_{\eta}$ takes one of the following normal forms:
\[
M_{\eta} \in \left\{ 
\begin{pmatrix}
\alpha & 0 \\ 
0 & \beta
\end{pmatrix}, 
\begin{pmatrix}
\alpha & 1 \\ 
0 & \alpha
\end{pmatrix}
\right\}.
\]
So if $\eta \in \big\{\eta_1, \eta_2, \eta_3, \eta_4, \eta_5\big\}$, we have:
\begin{align*}
M_{\eta_1} &= 
\begin{pmatrix} 
\alpha_1 & \frac{1 - \alpha_1^2}{\alpha_2} \\ 
\alpha_2 & -\alpha_1 
\end{pmatrix},
&
M_{\eta_2} &= 
\begin{pmatrix} 
1 & \alpha_1 \\ 
0 & -1 
\end{pmatrix},
&
M_{\eta_3} &=
\begin{pmatrix} 
-1 & \alpha_1 \\ 
0 & 1 
\end{pmatrix},
&
M_{\eta_4} &= 
\begin{pmatrix} 
1 & 0 \\ 
0 & 1 
\end{pmatrix},
&
M_{\eta_5} &= - M_{\eta_4}.
\end{align*}

From the above normal forms, the conditions $\operatorname{tr}(M_{\eta_i})=0$ and $\det(M_{\eta_i})=-1$ $\big($for $i=1,2,3\big)$ 
select exactly the normal form $\begin{pmatrix}1&0\\ 0&-1\end{pmatrix}$, and therefore the corresponding representative is precisely $\big(\Delta_{13},\ -\Delta_{23},\ 0\big)$. 
The matrix $M_{\eta_4}$ corresponds to the representative
$\big(\Delta_{13},\ \Delta_{23},\ 0\big)$, while 
$M_{\eta_5}$ corresponds to
$\big(-\Delta_{13},\ -\Delta_{23},\ 0\big)$. 
Since trace and determinant are invariants under similarity, 
these three representatives are not in the same orbit. 
Hence, we obtain the   algebras 
$\mathrm{B}_{16}^{1}$, $\mathrm{B}_{17}^{1}$, and $\mathrm{B}_{19}^{1}$.

\item
\underline{${\rm B}^{+}={\rm M}_{05}$.} Let 
$0\neq \theta =\big(
B_{1},B_{2},B_{3}\big) \in {\rm Z}^{2}\big( {\rm B
}^{+},{\rm B}^{+}\big) $. Then 
$\theta =\big( \alpha \Delta_{13},\ 0,\ 0\big) $ such that $\alpha ^{2}=1$. The automorphism group of $%
{\rm M}_{05}$  consists of the automorphisms $\phi $ given by a matrix of the following form:%
\begin{equation*}
\phi =%
\begin{pmatrix}
a_{11} & 0 & 0 \\ 
0 &1 & a_{23} \\ 
0 & 0 & 1%
\end{pmatrix}%
.
\end{equation*}%

Then $\theta \ast \phi =\theta $ for any $\phi \in \text{Aut}\big({\rm M}_{05}\big) $. So we get the algebras ${\rm B}_{20}$ and ${\rm B}_{21}$.

\item 
\underline{${\rm B}^{+}={\rm M}_{06}$.} Let 
$0\neq \theta =\big(
B_{1},B_{2},B_{3}\big) \in {\rm Z}^{2}\big( \rm{B%
}^{+},{\rm B}^{+}\big) $. Then $\theta =\big( \alpha \Delta
_{13},\alpha \Delta _{13}+\alpha \Delta _{23},0\big) $ such that $\alpha
^{2}=1$. The automorphism group of ${\rm M}_{06}$, $\text{Aut}\big( 
{\rm M}_{06}\big) $, consists of the automorphisms $\phi $ given by a
matrix of the following form:%
\begin{equation*}
\phi =%
\begin{pmatrix}
a_{11} & 0 & 0 \\ 
a_{21} & a_{11} & 0 \\ 
0 & 0 & 1%
\end{pmatrix}%
.
\end{equation*}%

Then $\theta \ast \phi =\theta $ for any $\phi \in \text{Aut}\big({\rm M}_{06}\big) $. So we get the algebras ${\rm B}_{22}$ and ${\rm B}_{23}$.

\item 
\underline{${\rm B}^{+}={\rm A}_{30}$.} Consider the change of basis 
$e_{1}^{\prime }=e_{1},$ 
$e_{2}^{\prime }=e_{2},$ 
$e_{3}^{\prime }={\rm i}e_{3}$. Then%
\begin{equation*}
{\rm A}_{30}\cong {\rm A}_{30}^{\prime }\ :\ 
e_{1}^{\prime}e_{1}^{\prime }=e_{1}^{\prime },\ 
e_{3}^{\prime }e_{3}^{\prime}=-e_{1}^{\prime }.
\end{equation*}%
So we assume ${\rm B}^{+}={\rm A}_{30}^{\prime }$. Let $0\neq\theta
=\big( B_{1},B_{2},B_{3}\big) \in {\rm Z}^{2}\big( 
{\rm B}^{+},{\rm B}^{+}\big) $. Then $\theta =\big( \alpha \Delta
_{13},\ 0,\ 0\big) $ such that $\alpha ^{2}=1$. The automorphism group of $%
{\rm A}_{30}^{\prime }$ consists of the automorphisms $\phi $ given by a matrix of the
following form:%
\begin{equation*}
\phi =%
\begin{pmatrix}
1 & 0 & 0 \\ 
0 & a_{22} & a_{23} \\ 
0 & 0 & \epsilon%
\end{pmatrix}%
:\ \epsilon ^{2}=1.
\end{equation*}%
Let $\phi =\bigl(a_{ij}\bigr)\in $ $\text{Aut}\big( {\rm A}
_{30}^{\prime }\big) $. Then 
$\theta \ast \phi =\big( \epsilon \alpha\Delta _{13},\ 0,\ 0\big) $. 
So we may choose $\epsilon $ such that $\epsilon
\alpha =1$ and hence we obtain the algebra ${\rm B}_{24}$.
\end{enumerate}


\section{The geometric classification of 
  algebras}

\subsection{Preliminaries: definitions and notation}
Given an $n$-dimensional vector space $\mathbb V$, the set ${\rm Hom}(\mathbb V \otimes \mathbb V,\mathbb V) \cong \mathbb V^* \otimes \mathbb V^* \otimes \mathbb V$ is a vector space of dimension $n^3$ and it has  the structure of the affine variety $\mathbb{C}^{n^3}.$ Indeed, let us fix a basis $e_1,\dots,e_n$ of $\mathbb V$. Then any $\mu\in {\rm Hom}(\mathbb V \otimes \mathbb V,\mathbb V)$ is determined by $n^3$ structure constants $c_{ij}^k\in\mathbb{C}$ such that
$\mu(e_i\otimes e_j)=\sum\limits_{k=1}^nc_{ij}^ke_k$. A subset of ${\rm Hom}(\mathbb V \otimes \mathbb V,\mathbb V)$ is {\it Zariski-closed} if it can be defined by a set of polynomial equations in the variables $c_{ij}^k$ $\big(1\le i,j,k\le n\big).$

Let $T$ be a set of polynomial identities.
The set of algebra structures on $\mathbb V$ satisfying polynomial identities from $T$ forms a Zariski-closed subset of the variety ${\rm Hom}(\mathbb V \otimes \mathbb V,\mathbb V)$. We denote this subset by $\mathbb{L}(T)$.
The general linear group ${\rm GL}(\mathbb V)$ acts on $\mathbb{L}(T)$ by conjugations:
$$ (g * \mu )(x\otimes y) = g\mu(g^{-1}x\otimes g^{-1}y)$$
for $x,y\in \mathbb V$, $\mu\in \mathbb{L}(T)\subset {\rm Hom}(\mathbb V \otimes\mathbb V, \mathbb V)$ and $g\in {\rm GL}(\mathbb V)$.
Thus, $\mathbb{L}(T)$ is decomposed into ${\rm GL}(\mathbb V)$-orbits that correspond to the isomorphism classes of algebras.
Let ${\mathcal O}(\mu)$ denote the orbit of $\mu\in\mathbb{L}(T)$ under the action of ${\rm GL}(\mathbb V)$ and $\overline{{\mathcal O}(\mu)}$ denote the Zariski closure of ${\mathcal O}(\mu)$.

Let $\bf A$ and $\bf B$ be two $n$-dimensional algebras satisfying the identities from $T$, and let $\mu,\lambda \in \mathbb{L}(T)$ represent $\bf A$ and $\bf B$, respectively.
We say that $\bf A$ degenerates to $\bf B$ and write $\bf A\to \bf B$ if $\lambda\in\overline{{\mathcal O}(\mu)}$.
Note that in this case we have $\overline{{\mathcal O}(\lambda)}\subset\overline{{\mathcal O}(\mu)}$. Hence, the definition of degeneration does not depend on the choice of $\mu$ and $\lambda$. If $\bf A\not\cong \bf B$, then the assertion $\bf A\to \bf B$ is called a {\it proper degeneration}. We write $\bf A\not\to \bf B$ if $\lambda\not\in\overline{{\mathcal O}(\mu)}$.

Let $\bf A$ be represented by $\mu\in\mathbb{L}(T)$. Then  $\bf A$ is  {\it rigid} in $\mathbb{L}(T)$ if ${\mathcal O}(\mu)$ is an open subset of $\mathbb{L}(T)$.
 Recall that a subset of a variety is called irreducible if it cannot be represented as a union of two non-trivial closed subsets.
 A maximal irreducible closed subset of a variety is called an {\it irreducible component}.
It is well known that any affine variety can be represented as a finite union of its irreducible components in a unique way.
The algebra $\bf A$ is rigid in $\mathbb{L}(T)$ if and only if $\overline{{\mathcal O}(\mu)}$ is an irreducible component of $\mathbb{L}(T)$.

\medskip

\noindent {\bf Method of the description of degenerations of algebras.} In the present work we use the methods applied to Lie algebras in \cite{GRH}.
First of all, if $\bf A\to \bf B$ and $\bf A\not\cong \bf B$, then $\mathfrak{Der}(\bf A)<\mathfrak{Der}(\bf B)$, where $\mathfrak{Der}(\bf A)$ is the   algebra of derivations of $\bf A$. We compute the dimensions of algebras of derivations and check the assertion $\bf A\to \bf B$ only for such $\bf A$ and $\bf B$ that $\mathfrak{Der}(\bf A)<\mathfrak{Der}(\bf B)$.

To prove degenerations, we construct families of matrices parametrized by $t$. Namely, let $\bf A$ and $\bf B$ be two algebras represented by the structures $\mu$ and $\lambda$ from $\mathbb{L}(T)$ respectively. Let $e_1,\dots, e_n$ be a basis of $\mathbb  V$ and $c_{ij}^k$ 
$\big(1\le i,j,k\le n\big)$ 
be the structure constants of $\lambda$ in this basis. If there exist $a_i^j(t)\in\mathbb{C}$ $\big(1\le i,j\le n$, $t\in\mathbb{C}^*\big)$ 
such that $E_i^t=\sum\limits_{j=1}^na_i^j(t)e_j$ ($1\le i\le n$) form a basis of $\mathbb V$ for any $t\in\mathbb{C}^*$, and the structure constants of $\mu$ in the basis $E_1^t,\dots, E_n^t$ are such rational functions $c_{ij}^k(t)\in\mathbb{C}[t]$ that $c_{ij}^k(0)=c_{ij}^k$, then $\bf A\to \bf B$.
In this case  $E_1^t,\dots, E_n^t$ is called a {\it parametrized basis} for $\bf A\to \bf B$.
In  case of  $E_1^t, E_2^t, \ldots, E_n^t$ is a {\it parametric basis} for ${\rm A}\to {\bf B},$ it will be denoted by
${\rm A}\xrightarrow{\big(E_1^t,\  E_2^t,\  \ldots, \ E_n^t\big)} {\bf B}$. 
To simplify our equations, we will use the notation $A_i=\langle e_i,\dots,e_n\rangle,\ i=1,\ldots,n$ and write simply $A_pA_q\subset A_r$ instead of $c_{ij}^k=0$ 
$\big(i\geq p$, $j\geq q$, $k< r\big).$


Let ${\rm A}(*):=\big\{ {\rm A}(\alpha) \big\}_{\alpha\in I}$ be a series of algebras, and let $\bf B$ be another algebra. Suppose that for $\alpha\in I$, $\bf A(\alpha)$ is represented by the structure $\mu(\alpha)\in\mathbb{L}(T)$ and $\bf B$ is represented by the structure $\lambda\in\mathbb{L}(T)$. Then we say that $\bf A(*)\to \bf B$ if $\lambda\in\overline{\big\{{\mathcal O}(\mu(\alpha))\big\}_{\alpha\in I}}$, and $\bf A(*)\not\to \bf B$ if $\lambda\not\in\overline{ \big\{{\mathcal O}(\mu(\alpha))\big\}_{\alpha\in I}}$.

Let $\bf A(*)$, $\bf B$, $\mu(\alpha)$ $\big(\alpha\in I\big)$ and $\lambda$ be as above. To prove $\bf A(*)\to \bf B$ it is enough to construct a family of pairs $\big(f(t), g(t)\big)$ parametrized by $t\in\mathbb{C}^*$, where $f(t)\in I$ and $g(t)\in {\rm GL}(\mathbb V)$. Namely, let $e_1,\dots, e_n$ be a basis of $\mathbb V$ and $c_{ij}^k$ $\big( 1\le i,j,k\le n \big)$ 
be the structure constants of $\lambda$ in this basis. If we construct $a_i^j:\mathbb{C}^*\to \mathbb{C}$ $\big(1\le i,j\le n\big)$ and $f: \mathbb{C}^* \to I$ such that $E_i^t=\sum\limits_{j=1}^na_i^j(t)e_j$ $\big(1\le i\le n\big)$ 
form a basis of $\mathbb V$ for any  $t\in\mathbb{C}^*$, and the structure constants of $\mu({f(t)})$ in the basis $E_1^t,\dots, E_n^t$ are such rational functions $c_{ij}^k(t)\in\mathbb{C}[t]$ that $c_{ij}^k(0)=c_{ij}^k$, then $\bf A(*)\to \bf B$. In this case  $E_1^t,\dots, E_n^t$ and $f(t)$ are called a parametrized basis and a {\it parametrized index} for $\bf A(*)\to \bf B$, respectively.

We now explain how to prove $\bf A(*)\not\to\mathcal  \bf B$.
Note that if $\mathfrak{Der} \ \bf A(\alpha)  > \mathfrak{Der} \  \bf B$ for all $\alpha\in I$ then $\bf A(*)\not\to\bf B$.
One can also use the following  Lemma, whose proof is the same as the proof of   \cite[Lemma 1.5]{GRH}.

\begin{lemma}\label{gmain}
Let $\mathfrak{B}$ be a Borel subgroup of ${\rm GL}(\mathbb V)$ and ${\rm R}\subset \mathbb{L}(T)$ be a $\mathfrak{B}$-stable closed subset.
If $\bf A(*) \to \bf B$ and for any $\alpha\in I$ the algebra $\bf A(\alpha)$ can be represented by a structure $\mu(\alpha)\in{\rm R}$, then there is $\lambda\in {\rm R}$ representing $\bf B$.
\end{lemma}

\subsection{The geometric classification of   
  algebras}

\subsubsection{Metabelian commutative  algebras}

\begin{theoremG1}\label{thm:geo-meta}
The variety of complex $3$-dimensional metabelian commutative algebras  
has dimension $8$ and it has $2$ irreducible components defined by  
\begin{center}
$\mathcal{C}_1=\overline{\mathcal{O}( {\rm M}_{07} )}$ \ and
$\mathcal{C}_2=\overline{\mathcal{O}( {\rm M}_{04}^{\alpha})}.$ 
\end{center}
In particular, there is only $1$ rigid algebra in this variety.
\end{theoremG1}

\begin{proof}
After carefully  checking  the dimensions of orbit closures of the more important for us algebras, we have 

\begin{longtable}{rclrclrclrcl}
$\dim \mathcal{O}({\rm M}_{04}^{\alpha})$&$=$&$8,$ &
$\dim \mathcal{O}({\rm M}_{07})$&$=$&$7.$  

\end{longtable}
\noindent

All necessary degenerations are given below 

\begin{longtable}{|lcr|}

\hline
    ${\rm M}_{02}  $&$ \xrightarrow{ \big(e_1+\frac{1}{2}e_2,\ e_3,\ te_2\big) } $&$ {\rm M}_{01}$    \\
\hline
   ${\rm M}_{03} $&$ \xrightarrow{ \big(te_1,\ e_2,\ te_3\big) } $&$ {\rm M}_{02}$    
\\
\hline

       ${\rm M}_{04}^{\alpha} $&$ \xrightarrow{ \big(e_1 + \frac{1}{2t\alpha}e_2 + te_3,\ 2te_1 + e_2,\ 2t^2(1-\alpha) e_1\big)} $&$  {\rm M}_{03}$    \\
\hline

          ${\rm M}_{04}^{\frac{t}{2}} $&$ \xrightarrow{ \big(e_1,\ te_1+e_2,\ \frac{1}{t}e_2+e_3\big) }$&$  {\rm M}_{05}$     
\\
\hline

     ${\rm M}_{04}^{1+t} $&$ \xrightarrow{ \big(e_1 +\frac{1}{t}e_2,\ e_2,\ e_3\big)} $&$  {\rm M}_{06}$ \\
\hline
\end{longtable}

${\rm M}_{04}^{\alpha} \not\to {\rm M}_{07}$ due to
$\mathcal R=\left\{
\begin{array}{lllllllll}
A_2^2=0 \\
\mbox{new basis for }{\rm M}_{04}^{\alpha}:
 \mathcal{B} =\big\{f_1=e_3, f_2=e_2, f_3=e_1\big\}
\end{array}\right\}.
$
\end{proof}

 \subsubsection{Derived commutative associative algebras}

\begin{theoremG2}\label{thm:geo-derassoccomm}
The variety of complex $3$-dimensional derived commutative associative algebras  
has dimension $12$ and it has $2$ irreducible components defined by 
\begin{center}
$\mathcal{C}_1=\overline{\mathcal{O}( {\rm J}_{07})}$ and
$\mathcal{C}_2=\overline{\mathcal{O}( {\rm A}_{01}^{\alpha, \beta, \gamma})}$.
\end{center}
In particular, there is only $1$ rigid algebra in this variety.
\end{theoremG2}

\begin{proof}
After carefully  checking  the dimensions of orbit closures of the more important for us algebras, we have 
\begin{longtable}{rclrclrclrcl}
$\dim\mathcal{O}({\rm A}_{01}^{\alpha, \beta, \gamma})$&$=$&$12,$  &
$\dim\mathcal{O}({\rm J}_{07})$&$=$&$9$. 
\end{longtable}

It is clear that no degeneration occurs from ${\rm A}_{01}^{\alpha, \beta, \gamma}$ to ${\rm J}_{07}$ because the dimension of the derived algebras cannot increase during degeneration.
All necessary degenerations are given below.

\begin{longtable}{|lcr|}
\hline
    ${\rm A}_{01}^{0,\ (1-\alpha){\rm i} {\Theta},\ {\alpha^2}{\Theta}^{-2}}  $&$ \xrightarrow{ \big(t e_1,\ te_1 + t e_2,\ e_1 + \alpha e_2 + {\rm i} {\Theta} e_3\big) } $&$ {\rm M}_{04}^{\alpha}$    \\
${\Gamma:=\sqrt{1-2 \alpha+2 \alpha^2}}$&&\\
\hline
    ${\rm A}_{01}^{0,\  \frac{4\rm i}{3},\  0}  $&$ \xrightarrow{ \big(t^2 e_1,\ 2te_1 + \frac{3t}{2}e_2,\ e_1 + \rm i e_3\big) } $&$ {\rm M}_{07}$    \\
\hline
    ${\rm A}_{01}^{t^{-1},\ \alpha t^{-1},\  0}  $&$ \xrightarrow{ \big(e_1,\ e_2,\ te_3\big) } $&$ {\rm A}_{02}^{\alpha}$    \\
\hline
    ${\rm A}_{05}^{0,\  \alpha}  $&$ \xrightarrow{ \big(e_1,\ \frac{1}{t}e_2,\ e_3\big) } $&$ {\rm A}_{03}^{\alpha}$    \\
\hline
    ${\rm A}_{05}^{t,\  \alpha}  $&$ \xrightarrow{ \big(e_1,\ (t^{-1}-\alpha)e_2+te_3,\ e_3\big) } $&$ {\rm A}_{04}^{\alpha}$    \\
\hline
 \multicolumn{2}{|l}{ ${\rm A}_{01}^{\alpha_{05},\  \beta_{05},\  \gamma_{05}} \  \xrightarrow{ \big(t^2 e_1+ (1-t^2) e_2,\ t e_1 + t e_2,\ \frac{t-1 - \alpha t + t^2 + \alpha t^3}{2 t^2-1} e_1 + \frac{t ( \alpha-1 + t - \alpha t^2)}{2 t^2-1} e_2+ \frac{\sqrt{2 t-1 - 2 \alpha t}}{2 t^2-1} e_3\big) } $}&$ {\rm A}_{05}^{\alpha, \beta}$    \\[1mm]
\multicolumn{3}{|l|}{
$\scriptstyle \alpha_{05}=\frac{t (t^2+t-1)-\alpha (1-t^2)^2}{t \sqrt{2 (1-\alpha) t-1}},$ \ 
$\scriptstyle \beta_{05}=\frac{t (\alpha t^2+t-1)}{\sqrt{2 (1-\alpha) t-1}},$ \
$\scriptstyle \gamma_{05}=\frac{2 \alpha^2 (t-3 t^3+2 t^5)+2 \alpha (1-t-2 t^2+2 t^3+3 t^4-2 t^6)+
t \big(2-4 t-2 t^2+4 t^3+\beta (2 t^2-1)^3\big)}{t +2 (\alpha-1) t^2}$}\\
    
\hline
    ${\rm A}_{01}^{t,\ t^{-3},\ \alpha}  $&$ \xrightarrow{ \big(e_1,\ t^2e_2,\ te_3\big) } $&$ {\rm A}_{06}^{\alpha}$    \\
\hline
    ${\rm A}_{01}^{1,\ \beta,\ \gamma}  $&$ \xrightarrow{ \big(e_1,\ te_2,\ te_3\big) } $&$ {\rm A}_{07}$    \\
\hline
    ${\rm A}_{01}^{\gamma t,\ \beta,\  \gamma}  $&$ \xrightarrow{ \big(e_1,\ t^2\gamma e_2,\ te_3\big) } $&$ {\rm A}_{08}$    \\
\hline
    ${\rm A}_{01}^{t,\ 0,\ \alpha t}  $&$ \xrightarrow{ \big(e_1,\ te_2,\ e_3\big) } $&$ {\rm A}_{09}^{\alpha}$    \\
\hline
    ${\rm A}_{01}^{0,\ t^{-3}\gamma^{-1},\ \gamma}  $&$ \xrightarrow{ \big(e_1,\ t^2\gamma e_2,\ te_3\big) } $&$ {\rm A}_{10}$    \\
\hline
    ${\rm A}_{01}^{0,\  \beta,\  t}  $&$ \xrightarrow{ \big(e_1,\ t e_2,\ e_3\big) } $&$ {\rm A}_{11}$    \\
\hline
  ${\rm A}_{01}^{\alpha_{12},\ \beta_{12},\ \gamma_{12}}  $&$ \xrightarrow{ \big( \frac{\alpha+t}{2 t}e_1+ \frac{t-\alpha}{2 t}e_2+\frac{\Theta}{2t}e_3 ,\ t (t^3-1)e_2+t+t^4e_1,\ (1+t^3)e_1\big) } $&$ {\rm A}_{12}^{\alpha, \beta}$    \\

\multicolumn{3}{|l|}{
$\scriptstyle \alpha_{12}=\frac{(\alpha+t) (1-t^3)}{(1+t^3){\Theta}},\  \beta_{12}= \frac{(\alpha-t) (1+t^3)}{(t^3-1){\Theta}}, \ 
\gamma_{12}=\frac{(t^3-1) \big(\alpha^2 (t^3-3)-2\alpha ( t-3 t^4)-t^2 (3+4 t-t^3-4t^7 -4 \beta(1- t^6))\big)}{(t^3+1) \big(\alpha^2 (t^3+3)-2 \alpha (t+3 t^4)+t^2 (3-4 t+t^3+4 t^7-4 \beta (1-t^6))\big)}$}\\

\multicolumn{3}{|l|}{
${\Theta} = 
\sqrt{\frac{\alpha^2 (3+t^3)-2 \alpha (t+3 t^4)+t^2 \big(3-4 t+t^3+4 t^7+4 \beta (t^6-1)\big)}{t^3-1}}$}\\

\hline
    ${\rm A}_{12}^{t,\  \alpha}  $&$ \xrightarrow{ \big(e_1,\ t^{-1} e_2,\ e_3\big) } $&$ {\rm A}_{13}^{\alpha}$    \\
\hline
    ${\rm A}_{12}^{0,\  \alpha}  $&$ \xrightarrow{ \big(e_1,\ t^{-1} e_2,\ e_3\big) } $&$ {\rm A}_{14}^{\alpha}$    \\
\hline
    ${\rm A}_{12}^{\alpha^{-\frac 12}t^{-\frac 32},\  \alpha^{-1} t^{-1}} $&$ \xrightarrow{ \big(e_1,\ t e_2,\ \alpha^{\frac 12} t^{\frac 12} e_3\big) } $&$ {\rm A}_{15}^{\alpha}$    \\
\hline
    ${\rm A}_{12}^{\alpha,\  t^{-2}} $&$ \xrightarrow{ \big(e_1,\ t^2 e_2,\ t e_3\big) } $&$ {\rm A}_{16}$    \\
\hline
    ${\rm A}_{12}^{\alpha,\  t^{-2}} $&$ \xrightarrow{ \big(e_1,\ e_2,\ t e_3\big) } $&$ {\rm A}_{17}$    \\
\hline
    ${\rm A}_{12}^{t^{-3},\  0} $&$ \xrightarrow{ \big(e_1,\ t^2 e_2,\ t e_3\big) } $&$ {\rm A}_{18}$    \\
\hline
    ${\rm A}_{12}^{\alpha,\  \beta} $&$ \xrightarrow{ \big(e_1,\ (\alpha t)^{-1} e_2,\ t e_3\big) } $&$ {\rm A}_{19}$    \\
\hline
    ${\rm A}_{12}^{0, \beta} $&$ \xrightarrow{ \big(e_1,\ t^2 e_2,\ t e_3\big) } $&$ {\rm A}_{20}$    \\
\hline
 \multicolumn{2}{|l}{ ${\rm A}_{01}^{
    (t-1){{\Theta}^{-1}}, \ 
    (t+1){{\Theta}^{-1}},\  
    (1+2 t-3 t^2+4\alpha t^4){\Theta}^{-2}} \  \xrightarrow{ \big(t e_1- t e_2,\ t^2 e_1 + t^2 e_2,\ \frac{t+1}{2t} e_1 + \frac{t-1}{2t} e_2+ \frac{{\Theta}}{2t} e_3\big) }$} &$ 
    {\rm A}_{21}^{\alpha}$    \\
${\Theta}:=\sqrt{1-2 t-3 t^2+4 \alpha t^4}$&&\\
\hline
    ${\rm A}_{21}^{0} $&$ \xrightarrow{ \big(t^{-1} e_1,\ t^{-2} e_2,\ t e_3\big) } $&$ {\rm A}_{22}$    \\
\hline
    ${\rm A}_{21}^{t^{-4}} $&$ \xrightarrow{ \big(t^{-1} e_1,\ t^{-2} e_2,\ t e_3\big) } $&$ {\rm A}_{23}$    \\
\hline

    ${\rm A}_{01}^{(\alpha-1){\Theta}^{-1}, \ 
    (\alpha-1){\Theta}^{-1},\  1} $&$ \xrightarrow{ \big(t e_1- t e_2,\ t^2 e_1 + t^2 e_2,\ \frac{\alpha+1}{2} e_1 + \frac{\alpha+1}{2} e_2+ \frac{\Theta}{2} e_3\big) } $&$ {\rm A}_{24}^{\alpha}$    \\
${{\Theta}:=\sqrt{1-2\alpha-3\alpha^2}}$&&\\
\hline
    ${\rm A}_{05}^{t^2,\  t^2} $&$ \xrightarrow{ \big(t e_1+t^{-1}e_2,\ -e_2,\ e_3\big) } $&$ {\rm A}_{25}$    \\
\hline
    ${\rm A}_{05}^{-\frac{t+1}4,\  0} $&$ \xrightarrow{ \big(t e_1-\frac{2t}{t+1} e_2,\ \frac{2t^2}{t+1}  e_2,\ 2e_3\big) } $&$ {\rm A}_{26}$    \\
\hline
    ${\rm A}_{12}^{-\beta,\  \beta} $&$ \xrightarrow{ \big(\beta t^2 e_1+t^2 e_2,\ \beta t^4  e_2,\ t e_3\big) } $&$ {\rm A}_{27}$    \\
\hline
    ${\rm A}_{12}^{t, -t} $&$ \xrightarrow{ \big(-t e_1+ e_2+ t e_3,\ -t  e_2,\ e_3\big) } $&$ {\rm A}_{28}$    \\
\hline
    ${\rm A}_{01}^{0,\ t^{-2},\  \gamma} $&$ \xrightarrow{ \big(e_1,\ t e_2,\ t e_3\big) } $&$ {\rm A}_{29}$    \\
\hline
    ${\rm A}_{01}^{0,\ \beta,\  0} $&$ \xrightarrow{ \big(e_1,\ t e_2,\ e_3\big) } $&$ {\rm A}_{30}$    \\
\hline

\end{longtable}

\end{proof}

\subsubsection{Derived Jordan algebras}

\begin{theoremG3}\label{thm:geo_derjor}
The variety of complex $3$-dimensional derived Jordan algebras has dimension $12$ and it has $7$ irreducible components defined by 
\begin{center}
$\mathcal{C}_1=\overline{\mathcal{O}( {\rm J}_{07})}$, \
$\mathcal{C}_2=\overline{\mathcal{O}( {\rm J}_{12})}$, \
$\mathcal{C}_3=\overline{\mathcal{O}( {\rm J}_{14})}$, \
$\mathcal{C}_4=\overline{\mathcal{O}( {\rm J}_{16})}$, \
$\mathcal{C}_5=\overline{\mathcal{O}( {\rm J}_{19})}$, \\
$\mathcal{C}_6=\overline{\mathcal{O}( {\mathcal J}_{01}^{\alpha})}$, \ and
$\mathcal{C}_7=\overline{\mathcal{O}( {\rm A}_{01}^{\alpha, \beta, \gamma})}$.
\end{center}
In particular, there are only $5$ rigid algebras in this variety.
\end{theoremG3}

\begin{proof}
After carefully  checking  the dimensions of orbit closures of the more important for us algebras, we have 
\begin{longtable}{rclrclrclrcl}
&&$\dim\mathcal{O}({\rm A}_{01}^{\alpha, \beta, \gamma})$&$=$&$12,$  & \\
&&$\dim\mathcal{O}({\mathcal J}_{01}^{\alpha})$&$=$&$10,$  & \\
&&$\dim\mathcal{O}({\rm J}_{07})$&$=$&$9$, \\
&&$\dim\mathcal{O}({\rm J}_{12})$&$=$&$8$, \\
$\dim\mathcal{O}({\rm J}_{16})$&$=$&$\dim\mathcal{O}({\rm J}_{19})$&$=$&$7$, \\
&&$\dim\mathcal{O}({\rm J}_{14})$&$=$&$3$.
\end{longtable}

Note that the derived algebras for ${\rm A}_{01}^{\alpha, \beta, \gamma}$ and ${\mathcal J}_{01}^{\alpha}$ are isomorphic to ${\mathfrak J}_{01}$ and ${\mathfrak J}_{04}$, respectively. But ${\mathfrak J}_{01}$ and ${\mathfrak J}_{04}$ are in different irreducible components in the variety of $2$-dimensional Jordan algebras. Since ${\mathfrak J}_{01} \not\rightarrow {\mathfrak J}_{04}$, we get ${\rm A}_{01}^{\alpha, \beta, \gamma} \not\rightarrow {\mathcal J}_{01}^{\alpha}$. Other non-degenerations can be obtained similarly to derived commutative associative algebras due to the dimension of the square of algebras:
$$ {\mathcal J}_{01}^{\alpha} \not\rightarrow 
\big\{ {\rm J}_{07},\ {\rm J}_{12},\ {\rm J}_{14},\ {\rm J}_{16},\ {\rm J}_{19} \big\}. $$ 

All necessary degenerations are given below 

\begin{longtable}{|lcr|}
\hline
    ${\mathcal J}_{01}^{(2{\rm i}\sqrt{2t^3})^{-1}}  $&$ \xrightarrow{ \big(e_1,\ \frac{{\rm i}}{\sqrt{2t}}e_2,\ \frac{{\rm i}t}{\sqrt{2t}} e_2 - {\rm i}\sqrt{2t} e_3\big) } $&$ {\mathcal J}_{02}$    \\
\hline
    ${\mathcal J}_{01}^{t^{-2}}  $&$ \xrightarrow{ \big(e_1,\ t e_1 + t^{-1}e_2,\ t e_3\big) } $&$ {\mathcal J}_{03}$    \\
\hline
    ${\mathcal J}_{01}^{\frac{2 (4+9\alpha)}{\sqrt{(9 \alpha-12)^3}}}  $&$ \xrightarrow{ \big(e_1 + \frac{2}{\sqrt{9 \alpha-12}} e_2,\ t e_1,\ \frac{2}{3} e_1 - \frac{4}{3 \sqrt{9 \alpha-12}}e_2 - \frac{\sqrt{9 \alpha-12}}{3} e_3\big) } $&$ {\mathcal J}_{04}^{\alpha}$    \\
\hline
    ${\mathcal J}_{01}^{\frac{-64}{\sqrt{(9t-48)^3}}}  $&$ \xrightarrow{ \big(e_1 + \frac{2}{\sqrt{9t-48}} e_2,\ t e_1,\ \frac{2}{3} e_1 - \frac{4}{3 \sqrt{9t-48}}e_2 - \frac{\sqrt{9t-48}}{3} e_3\big) } $&$ {\mathcal J}_{05}$    \\
\hline
    ${\mathcal J}_{01}^{t^{2}}  $&$ \xrightarrow{ \big(e_1,\ t e_1 + (t^3+t)e_2,\ \sqrt{t^2+1} e_3\big) } $&$ {\mathcal J}_{06}$    \\
\hline
    ${\mathcal J}_{01}^{t^{-2}}  $&$ \xrightarrow{ \big(e_1,\ e_2,\ t e_3\big) } $&$ {\mathcal J}_{07}$    \\
\hline
\end{longtable}

\end{proof}

\subsubsection{Bicommutative algebras}

\begin{theoremG4}\label{thm:geo_bicom}
The variety of complex $3$-dimensional bicommutative algebras has dimension $10$ and it has $4$ irreducible components defined by  
\begin{center}
$\mathcal{C}_1=\overline{\mathcal{O}( {\rm J}_{07})},$ \
$\mathcal{C}_2=\overline{\mathcal{O}( {\rm B}_{06}^{\gamma})},$ \
$\mathcal{C}_3=\overline{\mathcal{O}( {\rm B}_{07}^{\gamma})},$ \ and
$\mathcal{C}_4=\overline{\mathcal{O}( {\rm B}_{08}^{\gamma})}.$ 
\end{center}
In particular, there is only $1$ rigid algebra in this variety.
\end{theoremG4}

\begin{proof}

After carefully  checking  the dimensions of orbit closures of the more important for us algebras, we have 
\begin{longtable}{rcrcrcl} 
&&&&$\dim\mathcal{O}({\rm B}_{06}^{\gamma})$&$=$&$10,$ \\ 

$\dim\mathcal{O}({\rm J}_{07})$&$=$&$\dim\mathcal{O}({\rm B}_{07}^{\gamma})$&$=$&$\dim\mathcal{O}({\rm B}_{08}^{\gamma})
$&$=$&$9$.  
\end{longtable}

The algebras ${\rm B}_{07}^{\gamma}$ and ${\rm B}_{08}^{\gamma}$ are 
non-isomorphic, but opposite up to isomorphism; 
while ${\rm B}_{06}^{\gamma}$ is isomorphic to its opposite. 
All algebras obtained as degenerations from self-opposite algebras are also self-opposite, hence 
there are no degenerations from  ${\rm B}_{06}^{\gamma}$ to ${\rm B}_{07}^{\gamma}$ and ${\rm B}_{08}^{\gamma}.$
Thanks to~\cite{MS}, ${\rm J}_{07}$ is rigid in the variety of associative commutative algebras, and each commutative associative algebra is in the irreducible component defined by ${\rm J}_{07}$. Thus, it is obvious that ${\rm B}_{06}^{\gamma} \not\rightarrow \big\{ {\rm J}_{07},\ {\rm B}_{07}^{\gamma},\  {\rm B}_{08}^{\gamma} \big\}$.

All necessary degenerations are given below 

\begin{longtable}{|lcr|}
\hline
   ${\rm B}_{05}^{2t^{-1}-1}  $&$ \xrightarrow{ \big( -2 e_3,\  e_2,\  t e_1\big) } $&$ \mathcal{G}_{00}$    \\
\hline
    ${\rm B}_{06}^2  $&$ \xrightarrow{ \big(e_1+e_2,\ te_2,\   te_3\big) } $&$ {\rm B}_{01}$    \\
\hline
    ${\rm B}_{06}^2  $&$ \xrightarrow{ \big((t-t^2)e_1+te_2,\ t^3e_2+t^3e_3,\  t^2e_3\big) } $&$ {\rm B}_{02}$    \\
\hline
   ${\rm B}_{05}^{ 2t^{-1}-1}  $&$ \xrightarrow{ \big( t  e_1,\  t^2e_2,\  2t^2e_3\big) } $&$ {\rm B}_{03}$    \\
\hline
   ${\rm B}_{05}^\alpha  $&$ \xrightarrow{ \big( t e_1,\  e_2,\  te_3\big) } $&$ {\rm B}_{04}^\alpha$    \\
\hline
    ${\rm B}_{07}^{-2\alpha^{-1}-1}  $&$ \xrightarrow{ \big( \frac{2t+2t\alpha+t\alpha^2}{2(1+\alpha)}e_1 + (1+\alpha) e_2 + \frac{t\alpha^2}{2(1+\alpha)}e_3,\ t^2 e_1 + t\alpha^2 e_2,\   t^3 e_1\big) } $&$ {\rm B}_{05}^{\alpha}$ \\
\hline
    ${\rm B}_{06}^t  $&$ \xrightarrow{ \big(  e_1,\ t^2e_2,\  e_3\big) } $&$ {\rm B}_{09}$    \\
\hline
    ${\rm B}_{06}^{1+t}  $&$ \xrightarrow{ \big( e_1 + \frac{(1 + t)^2}{1 + 2 t} e_2,\ \frac{2 t (1 + t)^2}{1 + 2 t} e_2 + t e_3,\   e_3\big) } $&$ {\rm B}_{10}$ \\
\hline
    ${\rm B}_{06}^{(1+t^2)^{-1}}  $&$ \xrightarrow{ \big(  (1 + t^2)^2 e_1 + e_2 + t^2 (1 + t^2) e_3,\ -{t}(1+t^2)^{-1} e_2 + t e_3,\  (1 + t^2)e_3\big) } $&$ {\rm B}_{11}$ \\
\hline
    ${\rm B}_{14}  $&$ \xrightarrow{ \big(  t^{-1} e_1,\ t^{-2} e_2,\ e_3\big) } $&$ {\rm B}_{12}$ \\
\hline
    ${\rm B}_{15}  $&$ \xrightarrow{ \big(  t^{-1} e_1,\ t^{-2} e_2,\ e_3\big) } $&$ {\rm B}_{13}$ \\
\hline
    ${\rm B}_{07}^{1+t}  $&$ \xrightarrow{ \big( t e_1 + e_2,\ t^2 e_1,\ \frac{1+t}{2} e_1 + \frac{1}{2} e_3\big) } $&$ {\rm B}_{14}$ \\
\hline
    ${\rm B}_{08}^{-1-t}  $&$ \xrightarrow{ \big( t e_1 + e_2,\ t^2 e_1,\ \frac{1+t}{2} e_1 + \frac{1}{2} e_3\big) } $&$ {\rm B}_{15}$ \\
\hline
    ${\rm B}_{07}^\alpha  $&$ \xrightarrow{ \big(  e_2,\ t e_1,\ \frac{\alpha}{2}e_1 + \frac{1}{2}e_3\big) } $&$ {\rm B}_{16}^{\alpha}$ \\
\hline
    ${\rm B}_{07}^{-\alpha}  $&$ \xrightarrow{ \big( e_2,\ t e_1,\   \frac{\alpha}{2} e_1 + \frac{1}{2} e_3\big) } $&$ {\rm B}_{17}^{\alpha}$ \\
\hline
    ${\rm B}_{08}^{\alpha}  $&$ \xrightarrow{ \big( e_2,\ t e_1,\   \frac{\alpha}{2} e_1 + \frac{1}{2} e_3\big) } $&$ {\rm B}_{18}^{\alpha}$ \\ 
\hline
    ${\rm B}_{08}^{-\alpha}  $&$ \xrightarrow{ \big(  e_2,\ t e_1,\ \frac{\alpha}{2}e_1 + \frac{1}{2}e_3\big) } $&$ {\rm B}_{19}^{\alpha}$ \\
\hline
    ${\rm B}_{07}^t  $&$ \xrightarrow{ \big(  t e_2 +te_3,\ -t^2 e_1,\ e_3\big) } $&$ {\rm B}_{20}$ \\
\hline
    ${\rm B}_{08}^{-t}  $&$ \xrightarrow{ \big(  t e_2 +te_3,\ -t^2 e_1,\ e_3\big) } $&$ {\rm B}_{21}$ \\
\hline
    ${\rm B}_{14}  $&$ \xrightarrow{ \big(  t e_1,\ t e_2,\ e_3\big) } $&$ {\rm B}_{22}$ \\
\hline
    ${\rm B}_{15}  $&$ \xrightarrow{ \big(  t e_1,\ t e_2,\ e_3\big) } $&$ {\rm B}_{23}$ \\
\hline
    ${\rm B}_{06}^t  $&$ \xrightarrow{ \big(  (1 + t)^2 e_1 + t^2(1 + t)^2 e_2 + t (1 + t) e_3,\ t e_2,\  (1 + t)e_3\big) } $&$ {\rm B}_{24}$    \\
\hline

\end{longtable}

\end{proof}




\begin{thebibliography}{9}


 \bibitem{aak24}
 Abdelwahab H., 
 Abdurasulov K.,
 Kaygorodov I.,
The algebraic and geometric classification of noncommutative Jordan  algebras, Frontiers of Mathematics, 2025, DOI: 10.1007/s11464-024-0173-7


 
 \bibitem{afm}
 Abdelwahab H., 
 Fernández Ouaridi A., 
 Martín González C.,  
 Degenerations of Poisson algebras,  
Journal of Algebra and its Applications, 24 (2025),  3,  2550087. 


\bibitem{AKL}
 Abdelwahab H.,  Kaygorodov I.,  Lubkov R.,  
 The algebraic and geometric classification of $\delta$-Novikov algebras, 
 Boletín de la Sociedad Matemática Mexicana (3), 31 (2025),  3,  145.


\bibitem{akk}
 Abdurasulov K.,  Kaygorodov I.,     Khudoyberdiyev A.,  
 The algebraic classification of nilpotent  bicommutative    algebras,  
 Mathematics,  11  (2023), 3, 777.

 

\bibitem{kz}
 Abdurasulov K.,  Khudoyberdiyev A.,  Toshtemirova F.,
    The geometric classification of nilpotent  Lie-Yamaguti, Bol
  and compatible Lie algebras, 
    Communications in Mathematics,  33 (2025), 3, 10.



 \bibitem{jp2} 
 Alvarez M.,  Fehlberg Júnior R.,   Kaygorodov I.,  
 The algebraic and geometric classification of Zinbiel algebras,
 Journal of Pure and Applied Algebra, 226 (2022),  11,  107106.


 \bibitem{bcz} Bai Y.,   Chen Y.,   Zhang Z.,   Gelfand-Kirillov dimension of bicommutative algebras, Linear and Multilinear Algebra, 70 (2022),  22, 7623--7649.


 \bibitem{BBR} 
 Barreiro E.,   Benayadi S.,   Rizzo C., 
 Bialgebra theory for nearly associative algebras and ${\rm LR}$-algebras: equivalence, characterization, and ${\rm LR}$-Yang-Baxter equation, 
 Linear Algebra and its Applications, 711 (2025), 84--125.


 \bibitem{jp1} 
 Barrionuevo J., Tirao P.,   Sulca D.,  
 Deformations and rigidity in varieties of Lie algebras,  
 Journal of Pure and Applied Algebra, 227 (2023),   3,   107217.
 
 \bibitem{ben} Ben Hassine A.,  Chtioui T.,  Elhamdadi M.,  Mabrouk S.,  Cohomology and deformations of left-symmetric Rinehart algebras,
 Communications in Mathematics, 32 (2024),  2, 127--152.
 

 \bibitem{BO}
Benayadi S.,   Oubba H.,   
Nonassociative algebras of biderivation-type,  
Linear Algebra and its Applications, 701 (2024), 22--60. 
 

\bibitem{BBK}
Benayadi S.,  Boulmane S., Kaygorodov I., 
Nonassociative algebras of anti-biderivation-type, arXiv:2506.17228.


\bibitem{BG}
Bildanov R.,   Gorshkov I.,   
On $3$-generated axial algebras of Jordan type $\frac 12,$ 
Communications in Mathematics, 33 (2025),  3, 4.


\bibitem{BCZ}
Bokut L.,  Chen Y.,   Zhang Z.,   
Automorphisms of finitely generated free metabelian Novikov algebras, 
Journal of Algebra and its Applications, 25 (2026),  4,  2550370. 



\bibitem{BK}
  Burde D.,   Dekimpe K.,   Vercammen K.,  
  Complete ${\rm LR}$-structures on solvable Lie algebras, 
  Journal of Group Theory, 13 (2010),  5, 703--719. 


\bibitem{DKS} Burde D.,   Dekimpe K.,   Deschamps S., 
${\rm LR}$-algebras, 
New developments in Lie theory and geometry, 125--140, Contemp. Math., 491, Amer. Math. Soc., Providence, RI, 2009.


\bibitem{BC99} Burde D., Steinhoff C.,
    Classification of orbit closures of $4$--dimensional complex Lie algebras,
    Journal of Algebra, 214 (1999), 2, 729--739.



     \bibitem{DAS}
Dauletiyarova A.,   Abdukhalikov K.,  Sartayev B.,  
On the free metabelian Novikov and metabelian Lie-admissible algebras, 
Communications in Mathematics, 33 (2025),  3, 3.


\bibitem{DS}
 Dauletiyarova A.,  Sartayev B.,
Differential Novikov algebras, arXiv:2503.00457.

    
     \bibitem{DM}
De Medts T.,   Meulewaeter J.,   
Inner ideals and structurable algebras: Moufang sets, triangles and hexagons, 
Israel Journal of Mathematics, 259 (2024),  1, 33--88. 

     \bibitem{D19}
Drensky V.,  Varieties of bicommutative algebras, 
Serdica Mathematical Journal, 45 (2019),  2, 167--188.

     \bibitem{D23}
Drensky V.,   Invariant theory of free bicommutative algebras. Non-associative algebras and related topics, 231–243, Springer Proc. Math. Stat., 427, Springer, Cham, [2023], ©2023.

  
     \bibitem{dimz}
Drensky V.,   Ismailov N.,   Mustafa M.,   Zhakhayev B., 
Free bicommutative superalgebras, 
Journal of Algebra, 652 (2024), 158--187.


 \bibitem{DZ}
Drensky V.,   Zhakhayev B.,   Noetherianity and Specht problem for varieties of bicommutative algebras, 
Journal of Algebra, 499 (2018), 570--582.

   \bibitem{DI}
Dzhumadil'daev A.,   Ismailov N.,   
Polynomial identities of bicommutative algebras, Lie and Jordan elements, 
Communications in Algebra, 46 (2018),  12, 5241--5251.

 


 \bibitem{DIT}
Dzhumadildaev A., Ismailov N.,  Tulenbaev K.,  
Free bicommutative algebras, Serdica Mathematical Journal, 37 (2011),  1, 25--44. 

 \bibitem{DT}
Dzhumadildaev A., Tulenbaev K.,  
Bicommutative algebras,  
Russian Mathematical Surveys, 58 (2003),  6, 1196--1197.

  
 
 
  \bibitem{FKS}
Fehlberg J\'{u}nior R.,    Kaygorodov I.,  Saydaliyev A.,
The  geometric classification of symmetric Leibniz algebras, Communications in Mathematics, 33 (2025),  1, 10.

  \bibitem{FKS25}
Fehlberg J\'{u}nior R.,    Kaygorodov I.,  Saydaliyev A.,
The complete classification of irreducible components of varieties of Jordan superalgebras, Communications in Mathematics, 33 (2025),  3, 15.


\bibitem{jp3}
Fernández-Culma E.,   Rojas N.,   
On the classification of $3$-dimensional complex ${\rm Hom}$-Lie algebras, 
 Journal of Pure and Applied Algebra, 227 (2023),  5,  107272.



 
\bibitem{GRH}
Grunewald F.,  O'Halloran J.,
    Varieties of nilpotent Lie algebras of dimension less than six,
    Journal of Algebra, 112 (1988), 315--325.
 



 \bibitem{ims}
   Ismailov N.,  Mashurov F.,  Sartayev B., 
   On algebras embeddable into bicommutative algebras, 
Communications in Algebra, 52 (2024),  11, 4778--4785

  \bibitem{JVW}
  Jordan P., von Neumann J., Wigner E.,
On an algebraic generalization of the quantum mechanical formalism,
Annals of Mathematics (2), 35 (1934),  1, 29--64.



  \bibitem{jkk}
Jumaniyozov D.,   Kaygorodov I.,   Khudoyberdiyev A., 
The geometric classification of nilpotent commutative $\mathfrak{CD}$-algebras, 
Bollettino dell'Unione Matematica Italiana, 15 (2022),  3, 465--481. 

  \bibitem{k23}
  Kaygorodov I.,     
Non-associative algebraic structures: classification and structure, Communications in Mathematics,  32 (2024), 3, 1--62.



\bibitem{KMTZ}
Kaygorodov I.,   Mashurov F.,  Nam T., Zhang Z.,   
Products of commutator ideals of some Lie-admissible algebras, Acta Mathematica Sinica, English Series, 40 (2024),   8, 1875--1892.


\bibitem{MS}
Kaygorodov I.,  Khrypchenko M.,  Páez-Guillán P.,
The geometric classification of non-associative algebras: a survey, Communications in Mathematics, 32 (2024), 2, 185--284.
 

   \bibitem{kkp}
Kaygorodov I.,   Khrypchenko M.,   Popov Y.,  
The algebraic and geometric classification of nilpotent terminal algebras, 
Journal of Pure and Applied Algebra, 225 (2021),  6,   106625.

   \bibitem{kpv20}
Kaygorodov I.,   Páez-Guillán P.,   Voronin V., 
The algebraic and geometric classification of nilpotent bicommutative algebras, 
Algebras and Representation Theory, 23 (2020),  6, 2331--2347.

   \bibitem{kpv22}
 Kaygorodov I.,   Páez-Guillán P.,   Voronin V., 
 One-generated nilpotent bicommutative algebras, 
 Algebra Colloquium, 29 (2022),  3, 453--474.



   \bibitem{kg}
Kochloukova D.,   García Hernández  E.,  
Homological properties of metabelian restricted Lie algebras, 
Journal of Algebra, 643 (2024), 119--152. 


   \bibitem{ks25}
 Kolesnikov P.,   Sartayev B.,
White Manin product and Hadamard product,  arXiv:2503.21308v1

 

   \bibitem{l24}
 Lopes S.,
Noncommutative algebra and representation theory: symmetry, structure \& invariants, Communications in Mathematics,  32 (2024), 2, 63--117.
 
 \bibitem{ms}
Mashurov F.,   Sartayev B.,   
Metabelian Lie and perm algebras, 
Journal of Algebra and its Applications, 23 (2024),  4,  2450065.


   \bibitem{mv}
 Mishchenko S.,  Valenti A., 
Correspondence between some metabelian varieties and left nilpotent varieties,
Journal of Pure and Applied Algebra, 225 (2021),  3,  106538.

 
   \bibitem{OF} 
Öğüşlü N.,  Fındık Ş.,   
On invariants of free metabelian bicommutative algebras, 
Ricerche di Matematica, 74 (2025),  5, 3035--3050.


   \bibitem{p81}
Pchelintsev S., Varieties of algebras that are solvable of index 2, 
Mathematics of the USSR-Sbornik, 43 (1982), 2, 159--180.


   \bibitem{SS}
Siciliano S.,   Spinelli E.,  
Lie metabelian restricted universal enveloping algebras, 
Archiv der Mathematik, 84 (2005),  5, 398--405.

   \bibitem{SZ}
Shestakov I.,   Zhang Z.,  
Automorphisms of finitely generated relatively free bicommutative algebras, 
Journal of Pure and Applied Algebra, 225 (2021),  8,   106636. 



   \bibitem{T25}
Towers D.,   Nilpotency, solvability and Frattini theory for bicommutative, assosymmetric and Novikov algebras,  Communications in Algebra, 53 (2025),  3, 917--928.

\end{thebibliography}
\end{document}